\documentclass[final,onefignum,onetabnum]{siamonline190516}

\usepackage{booktabs}
\usepackage[mathscr]{euscript}
\usepackage{color,amsmath,amssymb,mathrsfs}
\usepackage{graphicx}
\usepackage{algorithmic,algorithm}
\usepackage{upgreek}
\usepackage{showkeys,cite}
\usepackage{multirow}
\usepackage{yfonts}
\usepackage{mathtools}
\usepackage[normalem]{ulem}
\usepackage{latexsym}
\DeclareMathAlphabet{\mathpzc}{OT1}{pzc}{m}{it}
\usepackage{url} 
\usepackage{comment}
\usepackage{dsfont}
\newcommand{\norm}[2][]{\left\Vert #2\right\Vert_{#1}}
\newcommand{\inprod}[2][]{\left\langle #2\right\rangle_{#1}}
\newcommand{\eps}{\varepsilon}
\newcommand{\bs}[1]{\ensuremath{\boldsymbol{#1}}}

\newtheorem{rmk}{Remark}
\newcommand{\afwd}{\mathcal{G}}
\newcommand{\cfwd}{\mathcal{F}}
\newcommand{\fwd}{\mathcal{R}}
\newcommand{\parm}{m}
\newcommand{\parb}{b} 
\newcommand{\both}{z}

\newcommand{\e}{{\bf e}}
\newcommand{\data}{{\bf d}}
\newcommand{\E}{\mathbb{E}}

\newcommand{\Der}[2][]{\mathsf{D}^{#1}_{#2}}

\newcommand{\G}{\Gamma}

\DeclareMathOperator{\diag}{diag}

\long\def\dontshow#1{{}}

\def\E{\mathbb{E}}

\DeclareMathOperator{\tr}{tr}

\begin{document}

\def\addressR{Optimization and Uncertainty Quantification, Sandia National Laboratories, Albuquerque, NM, 87185, USA}
\def\addressru{Department of Engineering Science, University of Auckland, Auckland 1010, New Zealand}
\def\addressruu{M\={a}ori and Pasifika Engineering Research Centre,  University of Auckland, Auckland 1010, New Zealand}
\def\addresswustl{Oden Institute for Computational Engineering \& Sciences, The University of Texas at Austin, Austin, TX, 78712, USA}
\def\addressucm{Applied Mathematics, School of Natural Sciences, University of California, Merced, CA, 95343, USA}

\author{Ruanui Nicholson\footnotemark[1] \footnotemark[2]
\and Radoslav Vuchkov\footnotemark[3]
  \and Umberto Villa\footnotemark[4]
  \and No\'{e}mi Petra\footnotemark[5]}
\renewcommand{\thefootnote}{\fnsymbol{footnote}}
\footnotetext[1]{\addressru\ (\email{ruanui.nicholson@auckland.ac.nz}).}

\footnotetext[2]{\addressruu}

\footnotetext[3]{\addressR\ (\email{rgvuchk@sandia.gov}).}
\footnotetext[4]{\addresswustl\ (\email{uvilla@austin.utexas.edu}).}
\footnotetext[5]{\addressucm\ (\email{npetra@ucmerced.edu}).}

\title{Taylor Approximation Variance Reduction for Approximation Errors\\ in PDE-constrained Bayesian Inverse Problems%
\thanks{Submitted to the editors \today.
    \funding{This research was supported by the Royal Society of New Zealand Te Ap\={a}rangi (Marsden Fund Council) grant MFP-24-UOA-279, by the U.S. National Science Foundation, Division of Mathematical Sciences grant CAREER-1654311, and by the U.S National Institutes of Health grants R01EB034261 and R01EB031585.}}}
      
\headers{Taylor Approximation Variance Reduction for Approximation Errors}{R.~Nicholson, R.~Vuchkov, U.~Villa, N.~Petra}

\maketitle
\begin{abstract}

In numerous applications, surrogate models are used as a replacement for accurate parameter-to-observable mappings when solving large-scale inverse problems governed by partial
differential equations (PDEs). The surrogate model may be a computationally cheaper alternative to the accurate parameter-to-observable mappings and/or may ignore additional unknowns or  sources of uncertainty. The Bayesian approximation error (BAE) approach provides a means to account for the induced uncertainties and approximation errors, i.e., the errors between the accurate parameter-to-observable mapping and the surrogate. The statistics of these errors are, however, in general unknown a priori, and are thus calculated using Monte Carlo sampling. Although the sampling is typically carried out offline, i.e., before considering the data, the process can still represent a computational bottleneck. In this work, we develop a scalable computational approach for reducing the costs associated with the sampling stage of the BAE approach. Specifically, we consider the Taylor expansion of the accurate and surrogate forward models with respect to the uncertain parameter fields either as a control variate for variance reduction or  as a means to directly and efficiently approximate the mean and covariance of the approximation errors. We propose efficient methods for evaluating the expressions for the mean and covariance of the Taylor approximations based on linear(-ized) PDE solves. Furthermore, the proposed approach is independent of the dimension of the uncertain parameter, depending instead on the intrinsic dimension of the data, ensuring scalability to high-dimensional
problems. The potential benefits of the proposed approach are demonstrated for two high-dimensional inverse problems governed by PDE examples, namely for  the estimation of a distributed Robin boundary coefficient in a linear diffusion problem, and for a coefficient estimation problem governed by a nonlinear diffusion problem.
\end{abstract}

\section{Introduction}
Estimation of parameters of interest based on indirect observed data is a common problem in numerous fields of science and engineering, and leads naturally to the study of inverse problems. The models linking the parameters of interest to the measurable data are often described using (sets of) partial differential equations (PDEs). Developing and solving these models can represent a significant challenge. The estimation of the parameters is often further complicated by the models containing additional unknown parameters or sources of uncertainty, such as material parameters, domain geometry, and/or boundary conditions. This leads to model errors that, if not properly accounted for, can lead to artifacts and biased estimate of the sought-after parameter. To this end, in the current paper, we assume the following formulation
\begin{equation}\label{eq: fwd0a}
\data=\afwd(\parm,\parb)+\e,
\end{equation}
where $\data$ is a vector of measured data, $\afwd$ represents a PDE-based model, $\parm$ denotes the parameters of primary interest, $\parb$ denotes additional parameters which are not of primary interest (i.e., {\em auxiliary}  or {\em nuisance} parameters), and $\e$ is a random vector represents measurement noise. We assume here that measurement noise and parameters are independent 
, though this assumption can be relaxed.

To expedite the inference procedure it is common to use a surrogate model in place of the accurate model linking the parameters of interest and the data. In this work, we use the term broadly to include reduced-order models (e.g., via projection-based model reduction), simplified physics models (e.g., setting uncertain parameters to nominal values), and data-driven emulators (e.g., Gaussian processes or neural networks)~\cite{asher2015review,castello2019modeling, frangos2010surrogate, zhu2018bayesian, heidenreich2015bayesian,babaniyi2021inferring,kolehmainen2011marginalization}. By substituting the accurate model $\afwd(\cdot,\cdot)$ for a surrogate model independent of $\parb$ (this could include setting $\parb$ to some nominal value), the relationship between the data and parameter can be written as
 \begin{equation}\label{eq: fwd0b}
\data\approx\cfwd(\parm)+\e,
\end{equation}
 where $\cfwd(\cdot)$ represents the surrogate model. However, replacing the accurate model with a surrogate model and/or fixing uncertain parameters to nominal values essentially always induces model errors and model uncertainties, which may lead to severely biased parameter estimates if they are not properly accounted for~\cite{kaipio2006statistical,KolehmainenBAE,kaipio2007statistical,nicholson2018estimation}. 

The Bayesian approach to inverse problems has become popular as it provides a framework to systematically account for, and quantify, uncertainty~\cite{kaipio2006statistical,stuart2010inverse,bardsley2018computational}. Within the Bayesian framework, the estimation problem is recast as problems of statistical inference, wherein any unknowns (such parameters, measurement locations, geometry, etc.) are treated as random variables, and the solution to the inverse problem is the posterior distribution, i.e., the conditional distribution of the unknowns given the measured data. The Bayesian approximation error (BAE) approach  is a particularly straightforward approach to incorporate modeling errors and uncertainties into the Bayesian framework~\cite{kaipio2013bayesian,kaipio2007statistical}. The BAE approach begins by rewriting (\ref{eq: fwd0a}) using the surrogate model $\cfwd$ but introducing an additional (additive) error term. Specifically, we write
 \begin{equation}\label{eq: fwd0c}
\data=\cfwd(\parm)+\bs{\eps}+\e,
\end{equation}
where $\bs{\eps}(\parm,\parb)=\afwd(\parm,\parb)-\cfwd(\parm)$, termed {\em the approximation error}, is used to account for the induced modeling errors and uncertainties. 

The BAE approach approximates the approximation errors $\bs{\eps}$ as conditionally Gaussian, i.e., $\bs{\eps}|\parm\sim\mathcal{N}(\bs{\eps}_{0|\parm},\Gamma_{\eps|\parm})$. The statistics of the approximation errors are not known a priori, and thus both the (conditional) mean $\bs{\eps}_{0|\parm}$ and (conditional) covariance $\Gamma_{\eps|\parm}$ are computed using Monte Carlo (MC) sampling \cite{KolehmainenBAE,kaipio2007statistical}. A key benefit to the BAE approach is that the MC sampling is carried out at the {\em offline} stage, i.e., prior to considering (or potentially collecting) the data, while at the inference stage only the surrogate model $\cfwd(\cdot)$ is evaluated, rather than the accurate model $\afwd(\cdot,\cdot)$. 
This is particularly critical in some applications where the number of MC samples required to accurately compute the second-order statistics of the approximation errors can be in the order of tens of thousands~\cite{tick2019modelling,kolehmainen2011marginalization}, or even hundreds of thousands~\cite{rimpilainen2019improved}. 
Furthermore, even for problems in which the number of required MC samples may be relatively low (in the order of hundreds or lower)~\cite{calvetti2015artificial,abdallh2012bayesian,nissinen2010compensation}, for particularly complex or large-scale models, only very few samples may be feasible to compute. 

To reduce the number of approximation error samples required we propose applying a control variate approach~\cite{robert1999monte,law2007simulation} to accurately compute the statistics of $\bs{\eps}$. 
Control variate-based approaches have been widely applied as a variance reduction technique for uncertainty quantification (in both forward problems and inverse problems), see for example~\cite{peherstorfer2018survey,dodwell2019multilevel}, and the references therein. In this work we investigate the use of both a linear and a quadratic Taylor approximation of $\bs{\eps}$ as a control variable to compute both $\bs{\eps}_{0|\parm}$ and $\bs\Gamma_{\eps|\parm}$. The same approach has been used to accelerate the solution of optimal control under uncertainty in~\cite{chen2019taylor}. However, in~\cite{chen2019taylor} the stochastic quantity of interest (QoI) was a scalar functional, while in the current paper the stochastic QoI is the approximation error which is, in general, multidimensional. An important feature of our approach is that the 
proposed control variate 
technique is agnostic to the specific construction of the surrogate model, requiring only that the surrogate be differentiable with respect to the parameters of interest. This flexibility allows the method to be integrated with a wide range of surrogate modeling strategies, making it suitable for diverse applications involving PDE-constrained inverse problems.

The remaining sections of this paper are organized as follows. We begin by providing in an overview of the Bayesian approach to inverse problems as well as a review of the BAE approach in Section~\ref{sec: BIP}. Section~\ref{sec: CV} covers the use of Taylor approximations as control variables for the computation of the mean and covariance of the approximation error. In Section~\ref{sec: Exs} we present two PDE-based numerical examples demonstrating the potential benefits of the proposed approach. Finally, Section~\ref{sec: Conc} provides concluding remarks.

\section{Bayesian Inverse Problems and Approximation Errors}\label{sec: BIP}
In this section, we present a summary of the Bayesian approach to infinite-dimensional PDE-constrained inverse problems and the Bayesian approximation error approach to account for model errors and uncertainties. We begin by introducing the notation used throughout the remainder of the paper. For a real Hilbert space $\mathscr{H}$, we denote the corresponding inner product by $\inprod[\mathscr{H}]{\cdot,\cdot}$ and the induced norm by $\norm[\mathscr{H}]{\cdot}$, i.e., for $v\in\mathscr{H}$ we have $\norm[\mathscr{H}]{v}=\inprod[\mathscr{H}]{v,v}^{1/2}$, while the topological dual is denoted $\mathscr{H}'$. For two real Hilbert spaces $\mathscr{H}_1$ and $\mathscr{H}_2$, we let $\mathscr{L}(\mathscr{H}_1,\mathscr{H}_2)$ denote the space of bounded linear operators from $\mathscr{H}_1$ to $\mathscr{H}_2$. Moreover, for a linear operator $\mathcal{A}\in\mathscr{L}(\mathscr{H}_1,\mathscr{H}_2)$, we denote by $\mathcal{A}^\ast\in\mathcal{L}(\mathscr{H}_2,\mathscr{H}_1)$ the adjoint. In the current paper we will be particularly interested in Gaussian measures on Hilbert spaces. As such, we use $\mathcal{N}(\parm_0,\mathcal{C}_{\parm})$ to denote a Gaussian measure with mean $\parm_0\in\mathscr{H}$ and covariance operator $\mathcal{C}_{\parm}$. Unless stated otherwise, we assume any infinite-dimensional covariance operators are strictly positive
self-adjoint trace-class operator, ensuring the Bayesian inverse problems considered are
well-defined~\cite{stuart2010inverse}. 

In finite dimensions, we denote a Gaussian distribution by $\mathcal{N}(\bs{\parm}_0,\bs{\G}_{\parm})$, with $\bs{\parm}_0\in\mathbb{R}^n$ the mean and $\bs{\G}_{\parm}\in\mathbb{R}^{n\times n}$ the symmetric positive definite covariance matrix. For $\bs{u}\in\mathbb{R}^n$ and $\bs{w}\in\mathbb{R}^m$, we define the outer product as $\bs{u}\otimes \bs{w}=\bs{u}\bs{w}^T$. The (cross-) covariance operator between $\bs{u}$ and $\bs{w}$ is then defined
\begin{equation}
\bs{\G}_{uw}=\mathbb{E}[(\bs{u}-\bs{u}_0)\otimes(\bs{w}-\bs{w}_0)]=\bs{\G}_{wu}^T.
\label{eq:crossCovs}
\end{equation}
Finally, when computing second order derivatives, we will make use of tensor contractions. Specifically, for $\bs{A},\bs{B}\in\mathbb{R}^{n\times n}$ we have $\bs{A}:\bs{B}=\tr(\bs{A}^T\bs{B})$, while for $\bs{\mathcal{T}}\in\mathbb{R}^{\ell\times n\times n}$ (which will represent a third order second derivative tensor) we define the contraction $\bs{\mathcal{T}}:\bs{A}\in\mathbb{R}^\ell$ component-wise as
\[[\bs{\mathcal{T}}:\bs{A}]_{i}=\tr(\bs{T}_i^T\bs{A})=\bs{T}_i:\bs{A},\]
for $i=1,2,\dots,\ell$, where each of the $\bs{T}_i=\bs{\mathcal{T}}[i,:,:]$ are the {\em slices} of $\bs{\mathcal{T}}$.

We consider the problem of estimating the parameter $\both=(\parm,\parb)\in\mathscr{H}$, with $\mathscr{H}=\mathscr{H}_\parm\times \mathscr{H}_\parb$, from data $\data\in\mathbb{R}^p$ with the set up
\begin{equation}\label{eq: fwd1}
\data=\afwd(\both)+\e,
\end{equation}
where $\afwd:\mathscr{H}\rightarrow \mathbb{R}^p$ denotes the parameter-to-observable (PtO) map and $\e\in\mathbb{R}^p$ denotes (random) noise in the data. We will be particularly interested in PDE-constrained inverse problems, and as such suppose (initially) that $\mathscr{H}_m=L^2(S)$ and $\mathscr{H}_\parb=L^2(T)$ with $S$ and $T$ bounded open sets of $\mathbb{R}^{d_i}$, with $d_i\in\{1,2,3\}$ and $i=1,2$. Moreover, we will assume that the PtO involves solving a PDE and then applying an observation operator. That is to say, the PtO map can be decomposed into the composition of two operators: 
\begin{equation}
\afwd(\cdot) := (\mathcal{B} \circ \mathcal{S})(\cdot),
\end{equation}
where $\mathcal{S}:\mathscr{H}\rightarrow \mathscr{W}$ denotes the PDE solution operator (for some suitably-chosen function space $\mathscr{W}$), and $\mathcal{B}:\mathscr{W}\to\mathbb{R}^p$ denotes the observation operator. It will be convenient to consider the underlying PDE in {\em residual form}. To this end, let $\mathcal V$ denote the test space and $\mathcal V'$ its topological dual, with duality pairing 
$\inprod[\mathcal{V}',\mathcal{V}]{\cdot,\cdot}=0$, and let $\mathcal{R}:\mathscr{W}\times\mathscr{H}\rightarrow\mathcal{V}'$ denote the residual operator. Then for a given $\both\in\mathscr{H}$, the state $u\in\mathscr{W}$ satisfies \begin{align}\label{eq: Rstrong}
   \fwd(u,\both)=0\quad \text{in }\mathcal{V}',
\end{align} 
which has the associated weak form
\begin{align}\label{eq: Rweak}
r(u,v,\both):=\inprod[\mathcal{V}',\mathcal{V}]{\fwd(u,\both),v}=0,\quad \forall v\in\mathcal{V}. 
\end{align}

Within the Bayesian framework, the solution to the estimation problem is the posterior law $\mu^{\rm po}_\both$ which is described using Bayes' law. Specifically, in infinite-dimensions, Bayes' law is given in terms of the Radon-Nikodym derivative of $\mu^{\rm po}_\both$ with respect to the prior law $\mu^{\rm pr}_\both$, 
\begin{equation}
\frac{d\mu^{\rm po}_\both}{d\mu^{\rm pr}_\both} \propto \pi(\data\vert \both),
\end{equation} 
where $\pi(\data\vert \both)$ is the likelihood~\cite{stuart2010inverse}.
Under the assumption of Gaussian noise model, i.e., $\pi(\e)=\mathcal{N}(\e_{0},\bs{\Gamma}_{e})$, the likelihood is given by~\cite{kaipio2013bayesian,kaipio2006statistical}
\begin{align}\label{eq:contLike}
\pi(\data|z)&=\pi_{e}(\data-\afwd(\both))\propto\exp{\left(-\frac{1}{2}\norm[\bs{\Gamma}_{e}^{-1}]{\data-\afwd(\both)-\e_{0}}^2\right)}.
\end{align}
In the current paper we also assume that $\parm$ and $\parb$ are a priori independent, and thus the prior measure on $\both$ can be represented using the product measure $\mu^{\rm pr}_\both=\mu^{\rm pr}_\parm\times \mu^{\rm pr}_\parb$.

\subsection{Discretized formulation of the Bayesian inverse problem}
In this paper we  discretise the parameters $\parm$ and $\parb$ using continuous Lagrange basis functions $\{\phi_i\}_{i=1}^n$ and $\{\psi_j\}_{j=1}^q$, that is, we set $\parm_h(x)=\sum_{i=1}^{n}\parm_i\phi_i(x)$ and $\parb_h(x)=\sum_{j=1}^{q}\parb_i\psi_j(x)$. Thus, the unknowns are the coefficient vectors $\bs{\parm}=(\parm_1,\parm_2,\dots,\parm_n)^T\in\mathbb{R}^n$ and $\bs{\parb}=(\parb_1,\parb_2,\dots,\parb_q)^T\in\mathbb{R}^q$, which we stack as $\bs{\both}=(\bs{\parm}^T,\bs{\parb}^T)^T\in \mathbb{R}^{n+q}$. We closely follow the discretization approach described in in~\cite{bui2013computational,petra2014computational}, where further details can be found. 

Bayes' theorem in finite dimensions is given by (see for example~\cite{kaipio2006statistical,calvetti2007introduction,calvetti2023bayesian})
\begin{align}
\pi(\bs{\both}|\data)\propto\pi(\data|\bs{\both})\pi(\bs{\both}),
\end{align}
where $\pi(\bs{\both})$ and $\pi(\bs{\both}|\data)$ are the finite-dimensional approximation of the prior and posterior measures $\mu^{\rm pr}_\both$ and $\mu^{\rm po}_\both$ respectively. Above, $\pi(\bs{\both}|\data)$ denotes the discrete counterpart of the likelihood in \eqref{eq:contLike} and reads
\begin{align}\label{eq:discLike}
\pi(\data|\bs{z})\propto\exp{\left(-\frac{1}{2}\norm[\bs{\Gamma}_{e}^{-1}]{\data-\afwd(\bs{\both})-\e_{0}}^2\right)},
\end{align}
where, with a slight abuse of notation, we use $\afwd(\bs{\both})$ to denote the PtO map evaluated at the finite element function
corresponding to the vector $\bs{\both}=(\bs{\parm}^T,\bs{\parb}^T)^T$. 


As alluded to, we are particularly interested in the case of Gaussian prior measures, which for the discretized parameters is written $\pi(\bs{\both})=\mathcal{N}(\bs{\both}_0,\bs{\G}_\both)$. As $\parm$ and $\parb$ are assumed to be independent a priori, the prior covariance matrix is block diagonal, i.e., $\bs{\Gamma}_\both=\text{diag}(\bs{\Gamma}_\parm,\bs{\Gamma}_\parb)$. 
The covariance matrices $\bs{\Gamma}_\parm$ and $\bs{\Gamma}_\parb$ are defined by using the inverse of (discretized) elliptic operators. Specifically, we have $\bs{\G}_\parm=(\bs{A}_\parm\bs{M}_\parm\bs{A})_\parm^{-1}$ and $\bs{\G}_\parb=(\bs{A}_\parb\bs{M}_\parb\bs{A})_\parb^{-1}$, where   
\begin{equation}
\begin{aligned}\label{eq: priorcovs}
[\bs{M}_\parm]_{ij}&=\!\int_S\phi_i\phi_j\;dx,\; [\bs{A}_\parm]_{ij}=\!\int_S \gamma_\parm\phi_i\phi_j+\kappa_\parm\nabla\phi_i\cdot\bs{\Theta}_\parm\nabla\phi_j\;dx,\; i,j\in\{1,2,\dots,n\},\\
[\bs{M}_\parb]_{k\ell}&=\!\int_T\psi_k\psi_\ell\;ds,\; [\bs{A}_\parb]_{k\ell}=\!\int_T \gamma_\parb\psi_k\psi_\ell+\kappa_\parb\nabla\psi_k\cdot\bs{\Theta}_\parb\nabla\psi_\ell\;ds,\; k,\ell\in\{1,2,\dots,q\},
\end{aligned}
\end{equation}
see~\cite{bui2013computational,petra2014computational} for more details. 

By introducing the block-diagonal mass matrix $\bs{M}=\text{diag}(\bs{M}_\parm,\bs{M}_\parb)$ and stiffness matrix~$\bs{A}=\text{diag}(\bs{A}_\parm,\bs{A}_\parb)$, the disretised prior can be written as
\begin{align}
    \pi(\bs{\both})=\pi(\bs{\parm})\pi(\bs{\parb}) & \propto \exp{\left(-\frac{1}{2}\norm[\bs{M}_\parm]{\bs{A}_\parm(\bs{\parm}-\bs{\parm}_{0})}^2\right)}\exp{\left(-\frac{1}{2}\norm[\bs{M}_\parb]{\bs{A}_\parb(\bs{\parb}-\bs{\parb}_{0})}^2\right)}\nonumber\\
    &=\exp{\left(-\frac{1}{2}\norm[\bs{M}]{\bs{A}(\bs{\both}-\bs{\both}_{0})}^2\right)}.\label{eq:discPrior}
\end{align}
Finally, combining (\ref{eq:discLike}) and (\ref{eq:discPrior}) results in the discrete posterior:
\begin{align}
\pi(\bs{\both}|\data)\propto\exp{\left(-\frac{1}{2}\norm[\bs{\Gamma}_{e}^{-1}]{\data-\afwd(\bs{\both})-\e_{0}}^2-\frac{1}{2}\norm[\bs{M}]{\bs{A}(\bs{\both}-\bs{\both}_{0})}^2\right)}.
\end{align}

\subsection{Solving the Bayesian inverse problem}
Fully characterizing the posterior distribution is typically intractable for large-scale inverse problems. A computationally feasible alternative is to compute a local Gaussian approximation\footnote{The proposed approach could be used within a sampling-based approach, such as MCMC, to characterise the full posterior} to the posterior~\cite{bui2013computational,nicholson2018estimation}. Specifically, we can take the approximation $\pi(\bs{\both}|\data)\approx\mathcal{N}(\bs{\both}_{\rm MAP},\bs{\Gamma}_{\both|\data})$, where $\bs{\both}_{\rm MAP}$ is the maximum a posteriori (MAP) estimate, and $\bs{\Gamma}_{\both|\data}$ is the approximate posterior covariance matrix found by linearizing the forward model $\afwd$ about the MAP estimate. The MAP estimate is defined
\begin{align}
   \bs{\both}_{\rm MAP}=\operatornamewithlimits{argmin}_{\bs{\both}\in\mathbb{R}^{n+q}} \mathcal{J}(\bs{\both}),\quad \mathcal{J}(\bs{\both}):=\frac{1}{2}\norm[\bs{\Gamma}_{e}^{-1}]{\data-\afwd(\bs{\both})-\e_{0}}^2+\frac{1}{2}\norm[\bs{M}]{\bs{A}(\bs{\both}-\bs{\both}_{0})}^2, 
\end{align}
while the approximate posterior covariance operator is given by
\begin{align}
   \bs{\Gamma}_{\both|\data}=(\bs{G}^\ast(\bs{\both}_{\rm MAP})\bs{\Gamma}_{e}^{-1}\bs{G}(\bs{\both}_{\rm MAP})+\bs{\Gamma}_\both^{-1})^{-1},
\end{align}
where $\bs{G}(z_{\rm MAP})$ is the (discretized) Fr\'echet derivative of $\afwd$ evaluated at $z_{\rm MAP}$.

In various applications not all unknown or uncertain parameters in the forward model are of interest, see e.g.,~\cite{nicholson2018estimation,babaniyi2021inferring,tarvainen2009approximation,candiani2021approximation} among many. As such, we recall that the set of unknown parameters can be decomposed as  $\bs{\both}=(\bs{\parm},\bs{\parb})$. We take $\bs{\parm}\in\mathbb{R}^n$ to denote the primary parameters of interest and $\bs{\parb}\in\mathbb{R}^q$ to denote the auxiliary (not of primary interest) parameters. In such situations the underlying goal is to compute the marginal posterior:
\begin{align}
    \pi(\bs{\parm}|\data)=\int_{\mathbb{R}^q}\pi(\bs{\parm},\bs{\parb}|\data)\,d\bs{\parb}.
\end{align} 
As pointed out above, computation of $\pi(\bs{\parm},\bs{\parb}|\data)$ is generally considered infeasible. Consequently, one could consider computing the marginal approximate (Gaussian) posterior:
\begin{align}
    \pi(\bs{\parm}|\data)\approx\int_{\mathbb{R}^q}\mathcal{N}(\bs{\both}_{\rm MAP},\bs{\Gamma}_{\both|\data})\,d\bs{\parb}=\mathcal{N}(\bs{\parm}_{\rm MAP},\bs{\Gamma}_{\parm|\data}),\quad \bs{\parm}_{\rm MAP}=\bs{P}\bs{\both}_{\rm MAP},\quad \bs{\Gamma}_{\parm|\data}=\bs{P}\bs{\Gamma}_{\both|\data}\bs{P}^T,\nonumber
\end{align} 
where $\bs{P}=\begin{bmatrix}\bs{I} &\bs{0}\end{bmatrix}\in\mathbb{R}^{n\times (n+q)}$ with $\bs{I}\in\mathbb{R}^{n\times n}$ the identity matrix. However, the approximate posterior $\mathcal{N}(\bs{\both}_{\rm MAP},\bs{\Gamma}_{\both|\data})$ still requires the computation of the joint MAP estimate $\bs{\both}_{\rm MAP}$ and the (approximate) joint posterior covariance operator $\bs{\Gamma}_{\both|\data}$. For large-dimensional problems, with expensive-to-evaluate PtO maps, this can be time-consuming or intractable. 

To avoid carrying out the joint inference and further reduce the cost of inference, we consider introducing a surrogate model $\cfwd:\mathscr{H}_\parm\to\mathbb{R}^p$ which is independent of $\parb$, and rewrite the forward problem (\ref{eq: fwd1}) in terms of the surrogate\footnote{where $\cfwd(\bs{\parm})$ denotes the surrogate model evaluated at the finite element function
corresponding to the vector $\bs{\parm}$.}:
\begin{align}\label{eq: fwd2}
    \data=\afwd(\bs{\parm},\bs{\parb})+(\cfwd(\bs{\parm})-\cfwd(\bs{\parm}))+\e=\cfwd(\bs{\parm})+\bs{\eps}+\e,
\end{align}
where $\bs{\eps} = \bs{\eps}(\bs{\parm},\bs{\parb}):=\afwd(\bs{\parm},\bs{\parb})-\cfwd(\bs{\parm})$ is termed the {\em approximation error} (or in some contexts {\em model discrepancy}).


In numerous works it is shown that neglecting the approximation errors $\bs{\bs{\eps}}$ and carrying out inference of $\bs{\parm}$ using the surrogate model $\cfwd$ can result in (significantly) misleading results, see for example~\cite{candiani2021approximation,tick2019modelling,babaniyi2021inferring,huttunen2018approximation,nicholson2023global,rimpilainen2019improved,huttunen2007approximation,castello2019modeling} among many. In the current paper we consider accounting for the approximation errors using the BAE approach. We note that other approaches, such as the so-called Kennedy-O'Hagan (KOH) method, based on the use of a Gaussian processes (GP), have also been applied to account for approximation errors~\cite{higdon2004combining,kennedy2001bayesian}. 

\subsection{The Bayesian approximation error approach}
In the BAE approach the approximation error is approximated as conditionally Gaussian, i.e., $
\bs{\eps}\vert\bs\parm\sim\pi(\bs{\eps}\vert\bs{\parm})\approx\mathcal{N}(\bs{\eps}_{0\vert\bs{\parm}},\bs{\Gamma}_{\eps\vert\parm})$, with conditional mean
\begin{align}
\bs{\eps}_{0\vert\parm}=\bs{\eps}_0+\bs{\Gamma}_{\eps\parm}\bs{\Gamma}_{\parm}^{-1}(\bs{\parm}-\bs{\parm}_0),\label{eq: epsMean}
\end{align}
and conditional covariance
\begin{align}
\bs{\Gamma}_{\eps\vert\parm}=\bs{\Gamma}_\eps -\bs{\Gamma}_{\eps\parm}\bs{\Gamma}_{\parm}^{-1}\bs{\Gamma}_{\parm\eps},\quad \bs{\Gamma}_{\eps\parm}=\mathbb{E}[(\bs{\eps}-\bs{\eps}_0)(\bs{\parm}-\bs{\parm}_0)^T]\bs{M}_\parm=\bs{\Gamma}_{\parm\eps}^\ast.\label{eq: epsCov}
\end{align} 
In many applications, the further approximation $\pi(\bs{\eps}\vert\bs{\parm})\approx\pi(\bs{\eps})=\mathcal{N}(\bs{\eps}_0,\bs{\G}_{\eps})$, referred to as the {\em enhanced error model}, is made~\cite{tarvainen2009approximation,babaniyi2021inferring,kaipio2013bayesian}. That is, the approximation errors are treated as independent of the parameter. This is generally a safe approximation in that the marginal covariance is always larger (in terms of quadratic forms) than the conditional, i.e., $\bs{\Gamma}_{\eps}\succeq\bs{\Gamma}_{\eps|\parm}$. Furthermore, this approximation can significantly reduce the computational costs of training the approximation errors, as discussed below.

The (statistics of the) approximation errors are generally not known a priori and must be calculated based on samples. Specifically, at the offline stage (i.e., prior to the collection or incorporation of data) the statistics of the approximation errors are computed using the MC approach on samples from the prior. That is, samples of the approximation errors are generated as 
\begin{align}
\bs{\eps}^{(\ell)}=\afwd(\bs{\parm}^{(\ell)},\bs{\parb}^{(\ell)})-\cfwd(\bs{\parm}^{(\ell)}),
\end{align} 
for $\ell=1,2,\dots N$, where $\bs{\both}^{(\ell)}= (\bs{\parm}^{(\ell)},\bs{\parb}^{(\ell)})$ are samples from the prior, i.e., $\bs{\both}^{(\ell)}\sim\pi(\bs{\both})$. The sample mean and covariance can then be computed:
\begin{align}
\bs{\eps}_{0}&=\frac{1}{N}\sum_{\ell=1}^N\bs{\eps}^{(\ell)},\quad
\bs{\Gamma}_\eps=\frac{1}{N-1}\sum_{\ell=1}^N \bs{E}\bs{E}^T\label{eq: BAEm},
\end{align} 
where $\bs{E}=[\bs{\eps}^{(1)}-\bs{\eps}_{0},\bs{\eps}^{(2)}-\bs{\eps}_{0},\dots,\bs{\eps}^{(N)}-\bs{\eps}_{0}]\in\mathbb{R}^{p\times N}$. 

\begin{rmk}
    
If the enhanced error model approximation is {\em not} made, the  (cross-)covariance matrix $\bs{\Gamma}_{\eps\parm}=\bs{\Gamma}_{\eps\parm}^\ast$ as well as the sample approximation of the prior covariance matrix must be computed, i.e., one also needs to calculate 
\begin{align}
\hat{\bs{\parm}}_{0}=\frac{1}{N}\sum_{\ell=1}^N\bs{\parm}^{(\ell)},\quad
\bs{\Gamma}_{\eps\parm}=\frac{1}{N-1}\sum_{\ell=1}^N \bs{E}\bs{S}^\ast=\bs{\Gamma}_{\parm\eps}^\ast,\quad \hat{\bs{\Gamma}}_\parm=\frac{1}{N-1}\sum_{\ell=1}^N \bs{S}\bs{S}^\ast,
\end{align} 
where $\bs{S}=[\bs{\parm}^{(1)}-\bs{\parm}_{0},\bs{\parm}^{(2)}-\bs{\parm}_{0},\dots,\bs{\parm}^{(N)}-\bs{\parm}_{0}]\in\mathbb{R}^{n\times N}$. As the inverse of $\hat{\bs{\Gamma}}_\parm$ is required\footnote{In general, use of the known prior covariance matrix rather ${\bs{\Gamma}}_\parm$ than the estimate $\hat{\bs{\Gamma}}_\parm$ leads to an indefinite conditional covariance matrix in~(\ref{eq: epsCov}), see e.g.,~\cite{kolehmainen2011marginalization}.} (see (\ref{eq: epsMean}) and (\ref{eq: epsCov})), for large-dimensional problems (accurate) computation of $\hat{\bs{\Gamma}}_\parm$ can become infeasible as typically it would require $N\geq n$.
\end{rmk}

With the mean and covariance of the (enhanced error model) approximation errors in hand, we can now
consider the inference of $\bs{\parm}$ using the forward problem (\ref{eq: fwd0a}) approximated by
\begin{align}\label{eq: fwd3}
    \data=\cfwd(\bs{\parm})+\bs{\eta}, \quad \bs{\eta}=\bs{\eps}+\e,\quad \bs{\eta}\sim\mathcal{N}(\bs{\eta}_0,\bs{\G}_\eta),\quad \bs{\eta}_0=\bs{\eps}_0+\e_0,\quad \bs{\G}_\eta=\bs{\G}_\eps+\bs{\G}_e,
\end{align}
where $\bs{\eta}$ is often referred to as the {\em total errors}, i.e, the sum of the noise and the approximation errors. Equation (\ref{eq: fwd3}) leads to an approximate marginal likelihood of the form
\begin{align}
\hat{\pi}(\data|\bs{\parm})=\pi_{\eta}(\data-\cfwd(\bs{\parm}))\propto\exp{\left(-\frac{1}{2}\norm[\bs{\Gamma}_{\eta}^{-1}]{\data-\cfwd(\bs{\parm})-\bs{\eta}_{0}}^2\right)},
\end{align}
and an approximate marginal posterior of the form 
\begin{align}
\hat{\pi}(\bs{\parm}|\data)\propto\exp{\left(-\frac{1}{2}\norm[\bs{\Gamma}_{\eta}^{-1}]{\data-\cfwd(\bs{\parm})-\bs{\eta}_{0}}^2-\frac{1}{2}\norm[\bs{M}_\parm]{\bs{A}_\parm(\bs{\parm}-\bs{\parm}_{0})}^2\right)}.
\end{align}
Furthermore, the Gaussian approximation to the approximate marginal posterior is then of the form $\hat{\pi}(\bs{\parm}|\data)\approx\mathcal{N}(\hat{\bs{\parm}}_{\rm MAP},\hat{\bs{\Gamma}}_{m|\data})$, with 
\begin{align}\label{eq:BAE_map}
   \hat{\bs{\parm}}_{\rm MAP}&=\operatornamewithlimits{argmin}_{\bs{\parm}\in\mathbb{R}^n}\frac{1}{2}\norm[\Gamma_{\eta}^{-1}]{\data-\cfwd(\bs{\parm})-\bs{\eta}_{0}}^2+\frac{1}{2}\norm[\bs{M}_\parm]{\bs{A}_\parm(\bs{\parm}-\bs{\parm}_{0})}^2,\\
   \hat{\bs{\Gamma}}_{\parm|\data}&=(\bs{F}^\ast(\hat{\bs{\parm}}_{\rm MAP})\bs{\Gamma}_{\eta}^{-1}\bs{F}(\hat{\bs{\parm}}_{\rm MAP})+\bs{\Gamma}_m^{-1})^{-1},\label{eq:BAE_cov} 
\end{align}
with $\bs{F}(\hat{\bs{\parm}}_{\rm MAP})$ the (discretized) Fr\'echet derivative of $\cfwd$ evaluated at $\hat{\bs{\parm}}_{\rm MAP}$, and $\bs{F}^\ast(\hat{\bs{\parm}}_{\rm MAP})$ the adjoint.

The mean and covariance of the approximation errors (as shown in~\ref{eq: BAEm}) can be computed before collecting  or considering the data, i.e., offline. However, for particularly large-scale and/or complex models the MC sampling procedure can in itself represent a significant computational undertaking. In an attempt to reduce the computational costs associated with computing the statistics of the approximation errors we consider the use of a control variate approach, which we describe next.

\section{Variance reduction of the approximation errors using Taylor approximations}\label{sec: CV}

We begin by recalling some key results on the use of the control variate approach as a variance reduction technique, and then develop the main approximation methods. More specifically, we consider the use of (first- and second-order) Taylor approximations of the approximation error as a control variable for computing the mean $\bs{\eps}_0$ and covariance $\bs{\G}_{\eps}$. The surrogate model $\cfwd$ is fairly arbitrary\footnote{Use of the Taylor expansions clearly requires some order of differentiability.}, as such, for simplicity of the presentation, we give details only for the Taylor approximations of the accurate PtO $\afwd$.

In what follows, we assume $\afwd:\mathscr{H}\to\mathbb{R}^p$ is twice continuously (Fr\'echet) differentiable. For a given point $\both_\ast\in\mathscr{H}$ having discrete representation $\bs\both_\ast\in\mathbb{R}^{n+q}$ we denote the discretized first derivative (the Jacobian matrix) at $\bs\both_\ast$ by $\bs{G}(\bs\both_\ast)\in\mathbb{R}^{p\times (n+q)}$ and the discretized bounded symmetric bilinear second derivative  (Hessian tensor)  evaluated at $\bs\both_\ast$ by  $\bs{\mathcal{H}}(\bs\both_\ast)\in\mathbb{R}^{p\times (n+q)\times (n+q)}$. The linear and quadratic Taylor series approximations are then given by
\begin{align}\label{eq: linBAE}
    \eps^{\rm lin}(\bs{\both})&=\afwd(\bs{\both}_\ast)+\bs{G}(\bs{\both}_\ast)(\bs{\both}-\bs{\both}_\ast),\\
    \eps^{\rm quad}(\bs{\both})&=\afwd(\bs{\both}_\ast)+\bs{G}(\bs{\both}_\ast)(\bs{\both}-\bs{\both}_\ast)+\frac{1}{2}\bs{\mathcal{H}}(\bs{\both}_\ast):((\bs{\both}-\bs{\both}_\ast)\otimes (\bs{\both}-\bs{\both}_\ast)).\label{eq: quadBAE}
\end{align}
Taking $\bs{\both}_\ast=\bs{\both}_0=\mathbb{E}[\bs{\both}]$, the mean and covariance of linear Taylor approximation are given by 
\begin{align}
     \bs{\eps}_0^{\rm lin}=\mathbb{E}[\bs{\eps}^{\rm lin} (\bs{\both})]=\afwd(\bs{\both}_0),\quad\bs{\Gamma}^{\rm lin}=\bs{G}(\bs{\both}_0)\bs{\Gamma}_{\both}\bs{G}(\bs{\both}_0)^\ast\label{eq: linstats}
\end{align}
respectively, while the mean and covariance of the quadratic approximation are given (see for example~\cite[Theorem 3.2d.4]{provost1992quadratic}) by 
\begin{align}
        \bs{\eps}_0^{\rm quad}=\mathbb{E}[\bs{\eps}^{\rm quad} (\bs{\both})]=\bs{\eps}_0^{\rm lin}+\frac{1}{2}\bs{\mathcal{H}}(\bs{\both}_0):\bs{\Gamma}_{\both},\quad\bs{\G}^{\rm quad}=\bs{\Gamma}^{\rm lin}+\bs{Q},\quad Q_{ij}=\frac{1}{2}\tr(\bs{H}_i\bs{\Gamma}_{\both}\bs{H}_j\bs{\Gamma}_{\both}),\label{eq: quadmean}
\end{align}
where $i,j\in\{1,2,\dots,p\}$, and the slices $\bs{H}_i=\bs{\mathcal{H}}(\bs{\both}_0)[i,:,:]\in\mathbb{R}^{(n+q)\times(n+q)}$ are symmetric matrices.

\subsection{Taylor approximations as control variables}

The control variate method is a well-established variance reduction technique for computing the expectation of a QoI ~\cite{robert1999monte,ross2022simulation,law2007simulation}. The success of the approach relies on having an additional random (control) variable with known expectation which is correlated with the QoI. The expectation of the QoI can then be computed based on MC samples of the difference between the QoI and the control while correcting using the known expectation of the control. In our case, the QoI is the approximation error $\bs{\eps}$, while the control variables are the first- and second-order Taylor approximations of $\bs{\eps}$, the details of which follow.

Using the linear approximation as a control variable, we can rewrite the mean of the approximation errors in \eqref{eq: BAEm} as
\begin{align}
\label{eq:eps0_linearCV}
\bs{\eps}_{0}=\E[\bs{\eps}-\bs{\eps}^{\rm lin}]+\bs{\eps}^{\rm lin}_0&\approx\frac{1}{N}\sum_{\ell=1}^N(\bs{\eps}^{(\ell)}-\bs{\eps}^{\rm lin}(\bs{\both}^{(\ell)}))+\bs{\eps}^{\rm lin}_0.
\end{align} 
On the other hand, using the quadratic Taylor approximation as the control variable, the mean of the approximation errors can be written as
\begin{align}
\label{eq:eps0_quadCV}
\bs{\eps}_{0}=\E[\bs{\eps}-\bs{\eps}^{\rm quad}]+\bs{\eps}^{\rm quad}_0&\approx\frac{1}{N}\sum_{\ell=1}^N(\bs{\eps}^{(\ell)}-\bs{\eps}^{\rm quad}(\bs{\both}^{(\ell)}))+\bs{\eps}^{\rm quad}_0.
\end{align} 
Similarly\footnote{as the covariance can be written in terms of expectations}, for the covariance of the approximation errors (see \ref{eq: BAEm}), we can write
\begin{align}
\label{eq:cov_linearCV}
\bs{\G}_\eps&\approx\frac{1}{N-1}\sum_{\ell=1}^N((\bs{\eps}^{(\ell)}-\bs{\eps}_{0})(\bs{\eps}^{(\ell)}-\bs{\eps}_{0})^T-(\bs{\eps}^{\rm lin}(\bs\both^{(\ell)})-\bs{\eps}^{\rm lin}_{0})(\bs{\eps}^{\rm lin}(\bs\both^{(\ell)})-\bs{\eps}^{\rm lin}_{0})^T)+\bs{\G}_{\bs{\eps}^{\rm lin}},
\end{align} 
using the linear Taylor approximation as the control variable, or, using the quadratic Taylor approximation, 
\begin{align}
\label{eq:cov_quadCV}
\bs{\G}_\eps&\approx\frac{1}{N-1}\sum_{\ell=1}^N((\bs{\eps}^{(\ell)}-\bs{\eps}_{0})(\bs{\eps}^{(\ell)}-\bs{\eps}_{0})^T-(\bs{\eps}^{\rm quad}(\bs\both^{(\ell)})-\bs{\eps}^{\rm quad}_{0})(\bs{\eps}^{\rm quad}(\bs\both^{(\ell)})-\bs{\eps}^{\rm quad}_{0})^T)+\bs{\G}_{\eps^{\rm quad}}.
\end{align} 
\begin{rmk} For completeness, we also report that the cross-covariance $\G_{\eps\both}$, is given by
\begin{align}
\bs{\G}_{\eps\both}&\approx\frac{1}{N-1}\sum_{\ell=1}^N((\bs{\eps}^{(\ell)}-\bs{\eps}_{0})(\bs{\both}^{(\ell)}-\bs{\both}_{0})^T-(\bs{\eps}^{\rm lin}(\bs\both^{(\ell)})-\bs{\eps}^{\rm lin}_{0})(\bs{\both}^{(\ell)}-\bs{\both}_{0})^T)+\bs{\G}_{\eps^{\rm lin}\both},\nonumber
\end{align}
using the linear Taylor approximation as the control variable, and
\begin{align}
\bs{\G}_{\eps\both}&\approx\frac{1}{N-1}\sum_{\ell=1}^N((\bs{\eps}^{(\ell)}-\bs{\eps}_{0})(\bs{\both}^{(\ell)}-\bs{\both}_{0})^T-(\bs{\eps}^{\rm quad}(\bs\both^{(\ell)})-\bs{\eps}^{\rm quad}_{0})(\bs{\both}^{(\ell)}-\bs{\both}_{0})^T)+\bs{\G}_{\eps^{\rm quad}\both},\nonumber
\end{align}
when the quadratic Taylor approximation is taken as the control variable, where 
\[\bs{\G}_{\eps^{\rm lin}\both}=\bs{\G}_{\eps^{\rm quad}\both}=\bs{G}\bs{\G}_{\both}.\]
\end{rmk} 

\subsection{Computation of derivatives}
The target applications of this paper are large-scale inverse problems governed by (possibly nonlinear) PDEs. In such settings,  computation and storage of the derivatives required in (\ref{eq: linBAE}) and (\ref{eq: quadBAE}), i.e., $\bs{G}$ and $\bs{\mathcal{H}}$, can be computationally prohibitive. However, during sampling it is only the action of these derivatives which is needed, both of which can be evaluated using linearized (sensitivity) solves (i.e., in a matrix free fashion). For a more general discussion on efficient higher order derivatives of PDE-based models we refer the reader to~\cite[Section 2]{maddison2019automated} or~\cite[Section 3]{alger2020tensor}. 

As pointed out above, the surrogate model is fairly arbitrary, and thus computations involved for the derivatives are specific to the surrogate. As such, here we focus our attention on efficient derivatives of the accurate PDE-based model $\afwd$. 
We provide the strong form of the required derivatives, while the weak forms can be derived by considering the weak form of the residual equations (\ref{eq: Rweak}).

We are interested in computing action of the derivatives on samples from the prior  $\bs{\both}_{\ell}\sim\pi(\bs{\both})$ for $\ell\in\{1,2,\dots,N\}$. To do this we require several derivatives involving the discretized version of residual form of the PDE $\fwd$ (see (\ref{eq: Rstrong})).  To this end, we let $s$ denote the number of state variables after discretization, we take $\bs{D}_u\in\mathbb{R}^{s\times s}$ to denote the discretized version  of $\Der{u}\fwd(\bs{u}_\ast,\bs{\both}_\ast)$, i.e., the derivative of $\fwd$ with respect to ${u}$ evaluated at $(\bs{u}_\ast,\bs{\both}_\ast)$, and $\bs{D}_\both\in\mathbb{R}^{s\times (n+q)}$ to denote the discretized version of $\Der{\both}\fwd(\bs{u}_\ast,\bs{\both}_\ast)$, where $\bs{u}_\ast \in \mathbb{R}^s$ solves the forward problem $\fwd(\bs{u}_\ast,\bs{\both}_\ast)=0$ for $\bs\both_\ast\in\mathbb{R}^{n+q}$, and we let $\bs{B}\in\mathbb{R}^{p\times s}$ denote the discrete representation of the observation operator $\mathcal{B}$.  The action of the Jacobian on $\tilde{\bs{\both}}_{\ell} := \bs{\both}_{\ell} - \bs{\both}_\ast$ can then be computed as
\begin{equation}
\label{eq:jac_action}
\bs{G} \tilde{\bs{\both}}_{\ell} :=-\bs {B}\bs{r}_{\ell},
\end{equation}
where $\bs{r}_{\ell}\in\mathbb{R}^s$ solves the sensitivity problem 
\begin{equation}
\label{eq:first_order_sensitivity}
\bs{D}_u \bs{r}_{\ell} =-\bs{D}_\both \tilde{\bs{\both}}_{\ell}. 
\end{equation}
Letting $\bs{\mathcal{D}}_{uu}\in\mathbb{R}^{s\times s\times s}$, $\bs{\mathcal{D}}_{u \both }\in\mathbb{R}^{s\times (n+q)\times s}$, and $\bs{\mathcal{D}}_{\both \both}\in\mathbb{R}^{s\times (n+q)\times (n+q)}$ denote the discretized version of the second derivatives, $\mathsf{D}_{u}^2\fwd(\bs{u},\bs{\both}_\ast)$, $\mathsf{D}_{u \both }\fwd(\bs{u},\bs{\both}_\ast)$, and $\mathsf{D}_{\both}^2\fwd(\bs{u},\bs{\both}_\ast)$, respectively, the action of the Hessian on $\tilde{\bs{\both}}_{\ell}$ is
\begin{equation}
\label{eq:hessian_action_sym}
   \bs{\mathcal{H}}:(\tilde{\bs{\both}}_{\ell}\otimes \tilde{\bs{\both}}_{\ell} ) := -\bs{B}\bs{s}_{\ell},
\end{equation}
where $\bs{s}_{\ell}\in\mathbb{R}^s$ solves the second-order sensitivity problem 
\begin{equation}
\label{eq:second_order_sensitivity}
\bs{D}_u \bs{s}_{\ell} = -\bs{\mathcal{D}}_{uu}:(\bs{r}_{\ell} \otimes \bs{r}_{\ell} )+2\bs{\mathcal{D}}_{u\both}:( \tilde{\bs{\both}}_{\ell}\otimes\bs{r}_{\ell})+\bs{\mathcal{D}}_{\both\both}:(\tilde{\bs{\both}}_{\ell}\otimes \tilde{\bs{\both}}_{\ell})
\end{equation}
and $\bs{r}_{\ell} \in \mathbb{R}^s$ is the solution of the first order sensitivity problem \eqref{eq:first_order_sensitivity}.

\paragraph{Computational costs} In Table~\ref{tab: solves} we give the number of solves (forward and sensitivities) required for each sample of the BAE when using the standard MC approach as well as using the linear and quadratic Taylor expansions as control variables. Importantly, these are independent of the dimensions of the parameter and of the data. For a linear PDE forward model, building the linear and quadratic Taylor approximation control variates requires solving one and two additional PDEs, respectively. However, both \eqref{eq:first_order_sensitivity} and \eqref{eq:second_order_sensitivity} share the same left hand side operator of the linear forward problem $\fwd(\cdot,\bs{\both}_\ast)=0$, allowing to amortize the finite element assembly and preconditioner construction. For a nonlinear forward PDE model requiring $N_{\rm fwd}$ Newton steps to converge, the additional relative cost of constructing the Taylor approximation is further reduced as the cost solving \eqref{eq:first_order_sensitivity} and \eqref{eq:second_order_sensitivity} corresponds to one and, respectively two, amortized Newton steps. This means for a target accuracy of the modeling errors statistics, use of the first and second order Taylor approximation as control variate result in computational savings if the variance of the estimators \eqref{eq:eps0_linearCV}-\eqref{eq:cov_linearCV} and \eqref{eq:eps0_quadCV}-\eqref{eq:cov_quadCV} is reduced by at least a factor of $(1+1/N_{\rm fwd})$ and $(1+2/N_{\rm fwd})$, respectively, compared to the na\"ive estimator \eqref{eq: BAEm}.
\begin{table}[!ht]
\centering
\caption{The computational cost for evaluating an approximation error sample using MC, and the control variate-enhanced MC approaches each in terms of the number of (accurate) forward, and sensitivity PDE solves. Here, $N_{\rm fwd}$ denotes the number of linearized PDE solves to compute the forward solution (note $N_{\rm fwd}=1$ for a linear PDE forward model).}

\label{tab: solves}
\begin{tabular}{ |c|c|c|c|c|} 
\hline
Method & Forward & Sensitivity & Total Linearized PDEs to sample $\bs\eps$\\
\hline
Monte Carlo & 1 & 0 &  $N_{\rm fwd}$ \\ 
Linear control variable & 1 & 1 & $N_{\rm fwd}+1$ \\ 
Quadratic control variable & 1 & 2 & $N_{\rm fwd}+2$ \\
\hline
\end{tabular}
\end{table}

With the aim of reducing the number of forward simulations, we present an efficient approach for approximating $\bs{\Gamma}_\eta^{-1}$ with controllable accuracy by
constructing a low rank approximation based on a generalized eigenproblem. The approach is similar to the frequently used method used to approximate posterior covariance matrices (and Hessian matrices) in large scale inverse problems, see e.g.~\cite{bui2013computational,villa2021hippylib} and the references therein. 

\paragraph{Convergence of the approximation errors} Recall, the total errors $\bs{\eta}$ (see (\ref{eq: fwd3})) are given by the sum $\bs{\eta}=\bs{\eps}+\e$, which has covariance matrix $\bs{\G}_\eta=\bs{\G}_\eps+\bs{\G}_e$. The covariance of the noise $\bs{\G}_e$ is assumed known, while constructing $\bs{\G}_\eps$ is costly as it requires PtO evaluations. We begin by considering the following generalized (symmetric) eigenproblem:
\begin{align}
    \bs{\G}_\eps\bs v_i=\lambda_i\bs{\G}_e\bs v_i,\quad \lambda_1\geq\lambda_2\geq\dots\geq\lambda_p.
\end{align}
When the generalized eigenvalues of decay rapidly we consider computing a low-rank approximation based on only the $r$ largest eigenvalues. Specifically, by setting $\bs{V}_r=[\bs{v}_1,\bs{v}_2,\dots,\bs{v}_r]\in\mathbb{R}^{p\times r}$ and $\bs{\Lambda}_r=\diag(\lambda_1,\lambda_2,\dots,\lambda_r)\in\mathbb{R}^{r\times r}$, and employing the Sherman-Woodbury-Morrison formula, the inverse of the total error covariance matrix can be decomposed as
\begin{align}\label{eq:Gdeco}
\bs{\G}_\eta^{-1}=(\bs{\G}_e+\bs{\G}_\eps)^{-1}=\bs{\G}_e^{-1}-\bs{V}_r\bs{D}_r\bs{V}_r^T+\mathcal{O}\left(\sum_{i=r+1}^p\frac{\lambda_i}{1+\lambda_i}\right),
\end{align}
where $\bs{D}_r=\diag(\lambda_1/(1+\lambda_1),\lambda_2/(1+\lambda_2),\dots,\lambda_r/(1+\lambda_r))\in\mathbb{R}^{r\times r}$.

The decomposition in (\ref{eq:Gdeco}) provides a  means to decide when enough approximation error samples have been computed. Specifically, when the (generalized) eigenvalues become small compared to 1, the truncation error term also becomes small and thus, when the eigenvalues greater than 1 have converged we likely have an adequate number of samples\footnote{In practice we use a truncation value of less than 1, for example 0.1.}.

\begin{rmk}\label{rem: noisepre}
In many cases the noise covariance is of the form $\bs{\G}_e=\delta^2_e \bs{I}$, i.e., the white noise case with $\delta_e$ the noise level in the data. In such cases we can compare the spectrum of $\bs{\G}_\eps$ to $\delta^2_e$ (in practice we compare to $0.1\times\delta_e^2$) rather than comparing the generalized eigenvalues to 1
\end{rmk}


\subsection{A sample-free approximation}\label{sec:sampFree}

In some settings, it may be feasible to entirely avoid sampling when computing the approximation error statistics. Specifically, when the approximation error $\bs{\eps}$ is sufficiently smooth in $\bs{\both}$, and the prior distribution on $\bs{\both}$ is Gaussian with covariance matrix $\bs{\Gamma}_{z}$, the statistics of the approximation error can be approximated directly using a Taylor expansion. That is, using a linear approximation we set $\bs{\eps}_0=\bs{\eps}^{\rm lin}_0$ and $\bs{\G}_\eps=\bs{\G}_{\eps^{\rm lin}}$ (see (\ref{eq: linstats})), or using a quadratic approximation we would use $\bs{\eps}_0=\bs{\eps}^{\rm quad}_0$ and $\bs{\G}_\eps=\bs{\G}_{\eps_{\rm quad}}$ (see (\ref{eq: quadmean})).   

This approach results in a Gaussian approximation error model \emph{without the need for any offline sampling}, and thus bypasses the computational cost associated with MC or control variate-based estimation of approximation error statistics. It is particularly appealing in large-scale PDE-constrained inverse problems, where even a modest number of high-fidelity forward solves may be prohibitive. The quality of the resulting approximation depends on how well the Taylor expansion captures the behavior of $\bs{\eps}$ under the prior. Nonetheless, the approach enables a fully deterministic, sampling-free  construction of the approximation error model, which can be attractive in large-scale settings. We demonstrate the effectiveness of this sample-free strategy in the Robin boundary inversion problem in Section~\ref{sec: ExRobin}, where we show the approximation can yield accurate posterior estimates.


\section{Numerical Examples}\label{sec: Exs}

To demonstrate the applicability of the proposed approach, we present two PDE-based model problems. Specifically, we consider an inverse Robin boundary coefficient problem motivated by corrosion detection as well as the problem of estimating a spatially distributed coefficient in a nonlinear (semilinear) diffusion problem motivated by photo-acoustic tomography (PAT). For the first of these examples, we investigate the use of both linear and a quadratic Taylor approximations as control variables for estimating the statistics of the approximation errors, while for the second we consider the use of a linear Taylor approximation as a control variable.

\subsection{Example 1: Linear diffusion problem with quadratic Taylor approximation control variable}\label{sec: ExRobin}
The first problem we consider is estimating the (spatially) distributed Robin coefficient
field on an inaccessible part of the domain boundary based on potential measurements taken on an accessible part of the boundary. The setup of the problem is motivated by corrosion detection (see~\cite{inglese1997inverse,nicholson2021joint,nicholson2018estimation} and the references therein), in which case the Robin coefficient
field gives an indication for the state (level of corrosion) of the inaccessible part of the domain, while the measurements represent electric potential or temperature. Specifically, for $\Omega$ a bounded open set in $\mathbb{R}^2$, the problem is to estimate the Robin coefficient field $\parm>0$ on $\Gamma_{\rm I}$ from pointwise measurements of the potential $u$ at points in ${\Omega}$, where the potential satisfies
\begin{subequations}\label{eq: Rob1}
\begin{align}
-\nabla\cdot(\parb\nabla u)&=0&&\text{in }\Omega\\
(\parb\nabla u)\cdot \bs{n}&=g&&\text{on }\Gamma_{\rm A}\\
(\parb\nabla u)\cdot \bs{n}+\parm u&=0&&\text{on }\Gamma_{\rm I}
\end{align}
\end{subequations}
with $\Gamma_{\rm I}\cap\Gamma_{\rm A}=\emptyset$, $\Gamma_{\rm I}\cup\Gamma_{\rm A}=\partial\Omega$, $\parb>\parb_{\rm min}>0$ denoting the conductivity, $g$ the applied flux, and ${\bs n}$ the outward facing unit normal.

In numerous settings, see~\cite{nicholson2018estimation} and the references therein, the conductivity $\parb$ is only approximately known (at best), but is not of interest, i.e., the conductivity is an auxiliary unknown. As shown in~\cite{nicholson2018estimation,alexanderian2022optimal}, simply neglecting the uncertainty in $\parb$ by setting it to a nominal value in (\ref{eq: Rob1}) leads to poor parameter and uncertainty estimates of $\parm$, while applying the BAE approach can effectively account for the additional model uncertainty. In line with these previous works,  we take the accurate forward model $\afwd$ as the composition of an observation operator and the solution operator of (\ref{eq: Rob1}). On the other hand, we take the surrogate model $\cfwd$  to be the same as the accurate model but with a fixed value for $\parb$, i.e., $\cfwd(\parm):=\afwd(\parm,\parb_\ast)$, where $\parb_\ast$ denotes a (possibly well justified) nominal value of $\parb$.

In this example we consider using both a linear Taylor approximation and a quadratic as control variables for computing the statistics of the approximation errors. We thus require the action of both the first and second order Fr\'echet derivatives of $\eps$. We set $\Omega=(0,1)\times(0,1/4)$ for this example and again solve both the accurate forward model and the surrogate model using a Galerkin finite element method
with piecewise linear continuous Lagrange basis functions. The regular triangular mesh of the domain consists of 750 nodes and 1372 elements which coincides with the mesh used to discretise the auxiliary parameters $\parb$, while $\parm$ is discretized on 50 nodes along the boundary $\Gamma_{\rm I}$, which coincide with the restriction of the mesh used to discretise $\Omega$ to $\Gamma_{\rm I}$. As such, for this example we have $n=50$ and $q=1372$. 

\begin{figure}[t!]
\includegraphics[height=0.27\textwidth]{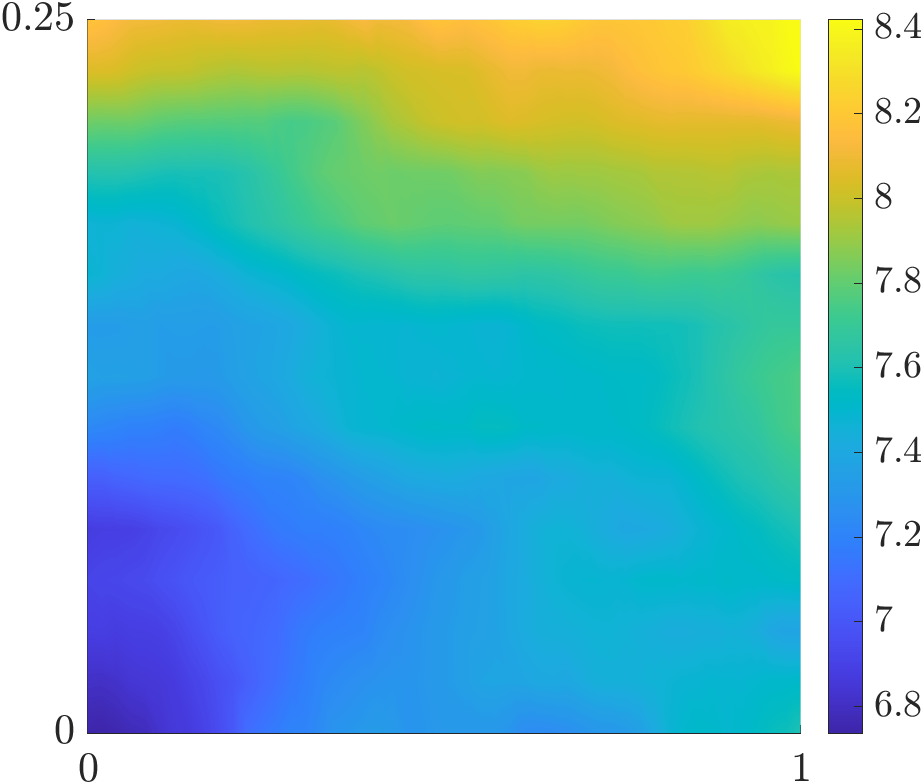} 
\hfill
\includegraphics[height=0.27\textwidth]{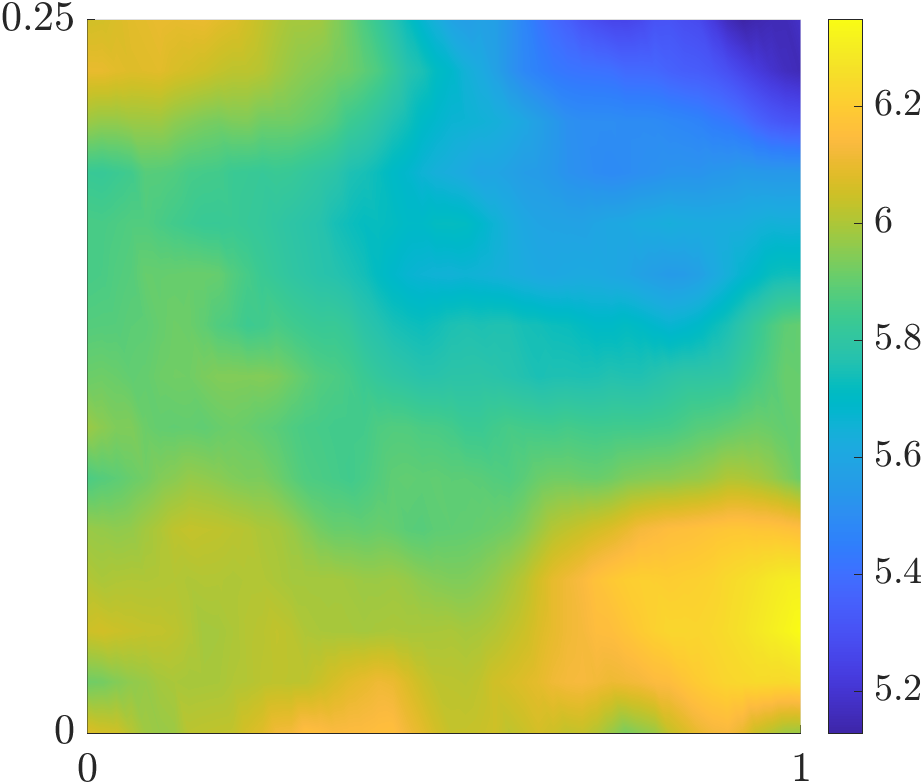}
\hfill
\includegraphics[height=0.27\textwidth]{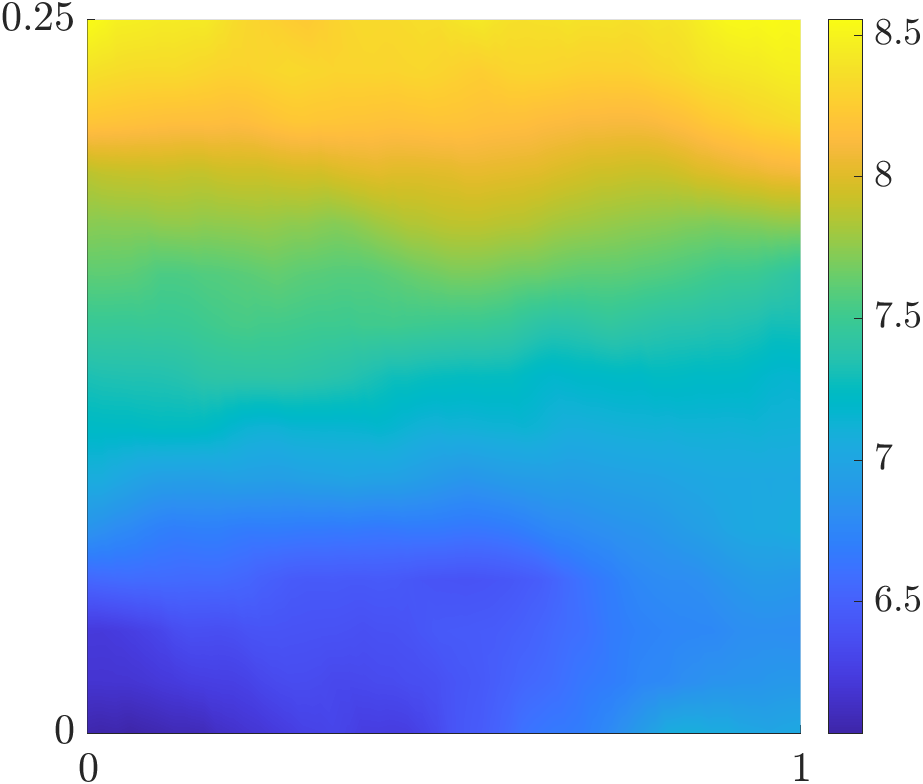}
    \caption{Example 1: Three prior samples of $\parb$.}
    \label{fig: priorSampsEx2}
\end{figure}

\begin{figure}[t!]
    \centering
\includegraphics[height=0.27\textwidth]{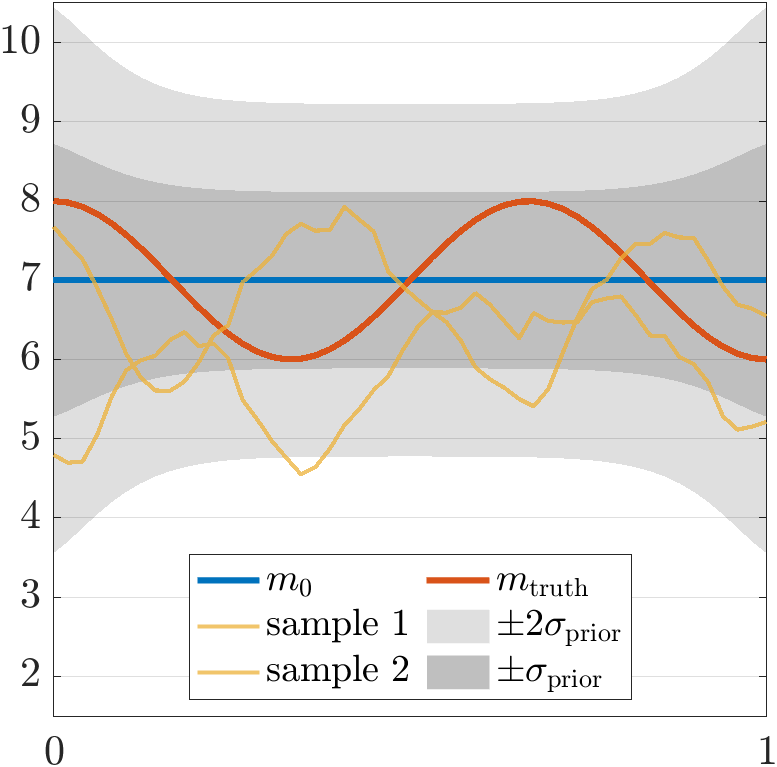}
\hfill
\includegraphics[height=0.27\textwidth]{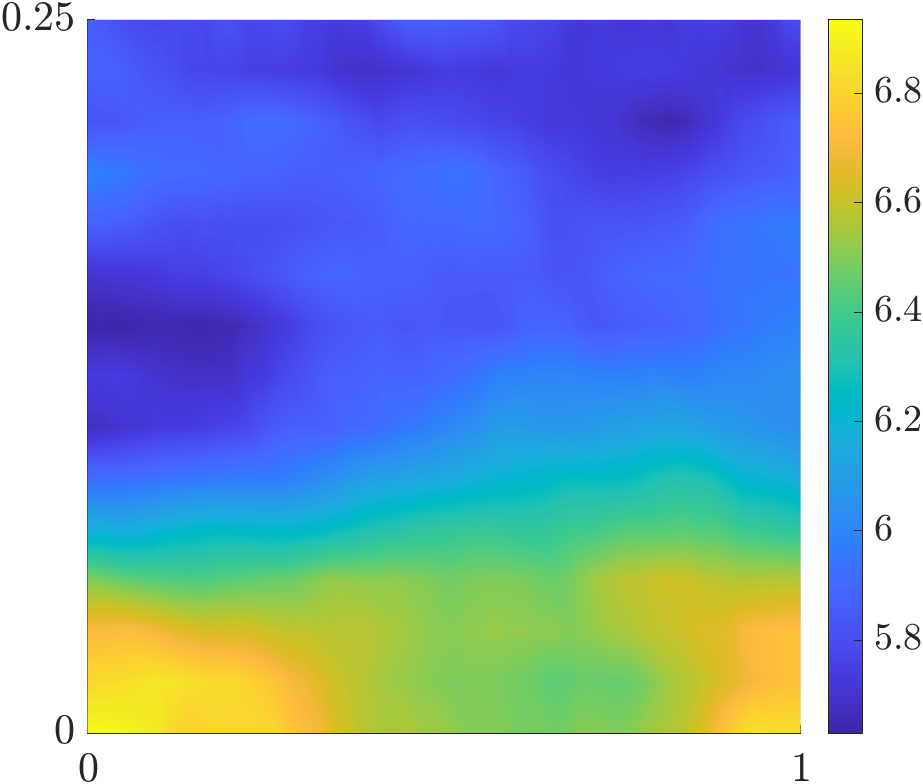}
\hfill
\includegraphics[height=0.27\textwidth]{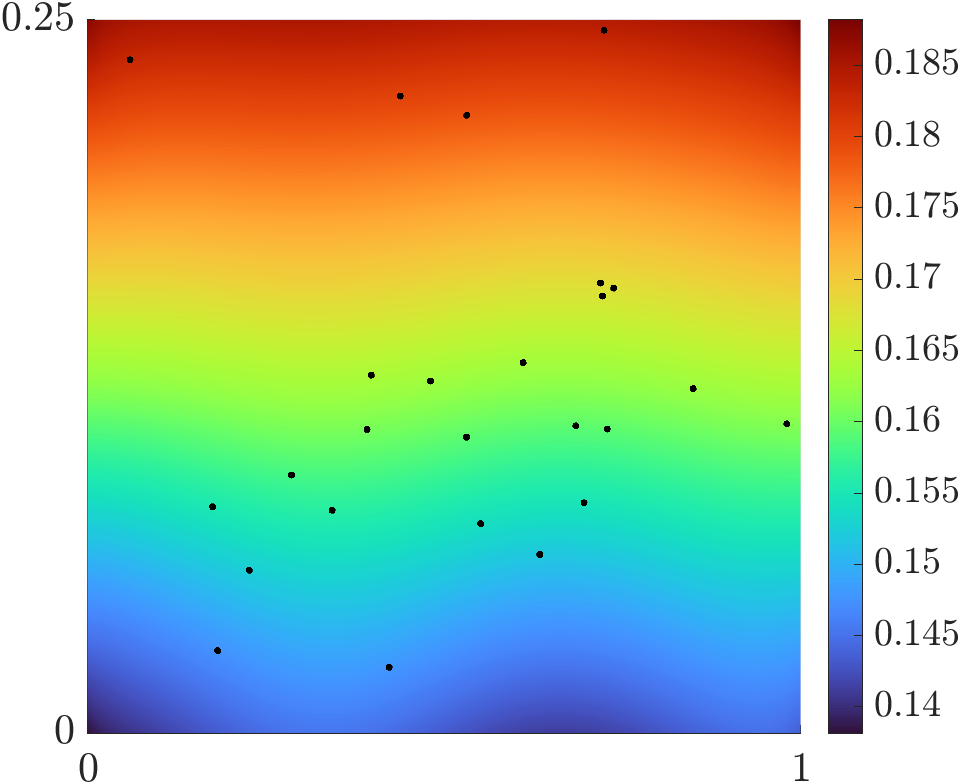}
    \caption{Example 1: The prior for the parameter $\parm_{\rm true}$ (left) with the true value $\parm_{\rm true}$ in red, two samples in yellow, the prior mean in blue and the $\pm\sigma_{\rm prior}$ and $\pm2\sigma_{\rm prior}$ intervals shaded. Also shown are the true auxiliary parameter $\parb_{\rm true}$ (centre) and the resulting solution to (\ref{eq: Rob1}) $u_{\rm true}$ (right) with measurement locations shown as black dots.}
    \label{fig: truthsEx2}
\end{figure}

The prior means for the uncertain parameters are set as $\parm_0=7$ and $\parb_0=7$, while the covariance operators (definied in~\ref{eq: priorcovs}) are set by taking $\gamma_\parm=10$, $\gamma_\parb=100$, $\kappa_\parm=0.1$, $\kappa_\parb=100$, $\Theta_\parm=1$ and $\bs{\Theta}_\parb=\diag(1,0.025)$. Three prior samples of $\parb$ are shown in Figure~\ref{fig: priorSampsEx2} while three samples of $\parm$ are shown in Figure~\ref{fig: truthsEx2} (left). We take $\parb_{\rm min}=10^{-6}$ and enforce the constraint $\parb>\parb_{\rm min}$ by discarding (and resampling) any prior samples which do not satisfy the condition. However, due to the choice of the parameters $\gamma_\parb$, $\kappa_\parb$, and $\bs{\Theta}_\parb$, all prior samples of $\parb$ satisfied the condition.

The inaccessible part of the boundary (where $\parm$ is defined) is defined $\Gamma_{\rm I}:=[0,1]\times \{0\}=S$ and we set the flux as with $g=1$. Lastly, data consists of noisy pointwise measurements of $u$ at 25 randomly chosen points throughout the domain, see Figure~\ref{fig: truthsEx2} with the noise being of the form $\e\sim\mathcal{N}(\bs{0},\delta_e^2\bs{I})$, where $\delta_e=(\max(\mathcal{B}u)-\min(\mathcal{B}u))\times 1/100$, i.e,  the noise level is set at 1\% of
the range of the noiseless measurements. 
 The true parameter of interest $\parm_{\rm true}$, as well as the true auxiliary parameter $\parb_{\rm true}$ along with the resulting solution $u$ from which data is extracted for inversions are shown in Figure \ref{fig: truthsEx2}.

The discussion on the convergence of the approximation error sampling as well as the inversion results for the numerical results is left to Section~\ref{sec: AllRes}.

\subsection{Example 2: Nonlinear diffusion problem with linear Taylor approximation control variable}
For the second problem we consider estimating $\parm$ within a (bounded open) domain $\Omega$ based on noisy measurements of $u$ which satisfies the semi-linear diffusion problem 
\begin{subequations}\label{eq: PAT1}
\begin{align}
-\Delta u+\exp(\parm)u+\exp{(\parb)}|u|u&=0\quad \text{in }\Omega\\       
         u&=g\quad \text{on }\Gamma_{\rm D},\\
         \nabla u\cdot {\bs n} &=0\quad \text{on }\Gamma_{\rm N},
\end{align}
\end{subequations}
with $\Gamma_{\rm D}\cap\Gamma_{\rm N}=\emptyset$ and $\Gamma_{\rm D}\cup\Gamma_{\rm N}=\partial\Omega$ and ${\bs n}$ the outward facing unit normal. Such a problem arises in, for example, two-photon photo-acoustic tomography~\cite{ren2018nonlinear,bardsley2018quantitative}, where $\parm$ represents the (log-) single-photon absorption coefficient, $\parb$ is the (log-) intrinsic two-photon absorption coefficient, and $g$ models the incoming near infra-red (NIR) photon source. 

We consider the problem of estimating  $\parm$ based on measurements of $u$, while $\parb$ is assumed to be not of interest but is unknown, i.e., is an auxiliary unknown. For any $\parb\in\mathscr{H}_\parb$ problem (\ref{eq: PAT1}) is non-linear, and typically requires iterative methods to solve for $u$. On the other hand, 
neglecting the nonlinear absorption term (i.e., taking $\exp(\parb)\rightarrow 0$) leads to a linear elliptic PDE which can be solved directly, and provides the motivation for the particular surrogate model in this example. More precisely, we take the accurate PtO $\afwd(\parm,\parb)$ to represent solving (\ref{eq: PAT1}) for $u$ and applying the observation operator $\mathcal{B}$ to the solution. On the other hand, the surrogate model $\cfwd(\parm)$ represents solving
\begin{subequations}\label{eq: PAT2}
\begin{align}
-\Delta u+\exp(\parm)u&=0\quad \text{in }\Omega\\       
         u&=g\quad \text{on }\Gamma_{\rm D},\\
         \nabla u\cdot {\bs n} &=0\quad \text{on }\Gamma_{\rm N},
\end{align}
\end{subequations}
that is, neglecting the nonlinear absorption in $(\ref{eq: PAT1})$, for $u$  and applying the (same) observation operator. In this example we only consider using a linear Taylor approximation (see Equation~(\ref{eq: linBAE})) as a control variable.

\begin{figure}[t!]
$\parm$ \hrulefill\vspace{5pt}

\includegraphics[height=0.27\textwidth]{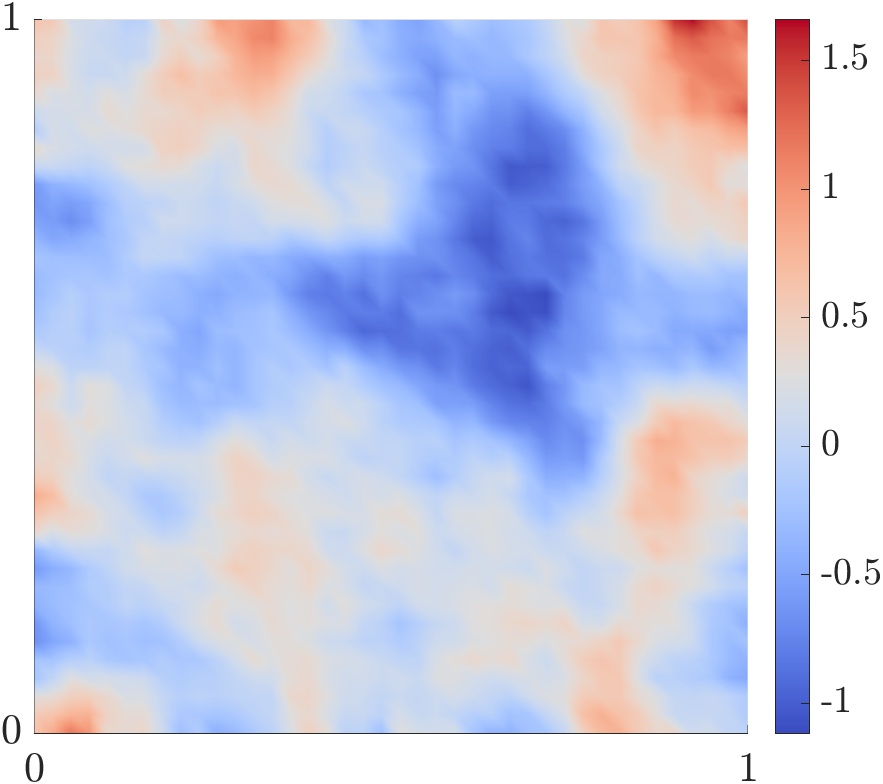}   
\hfill
\includegraphics[height=0.27\textwidth]{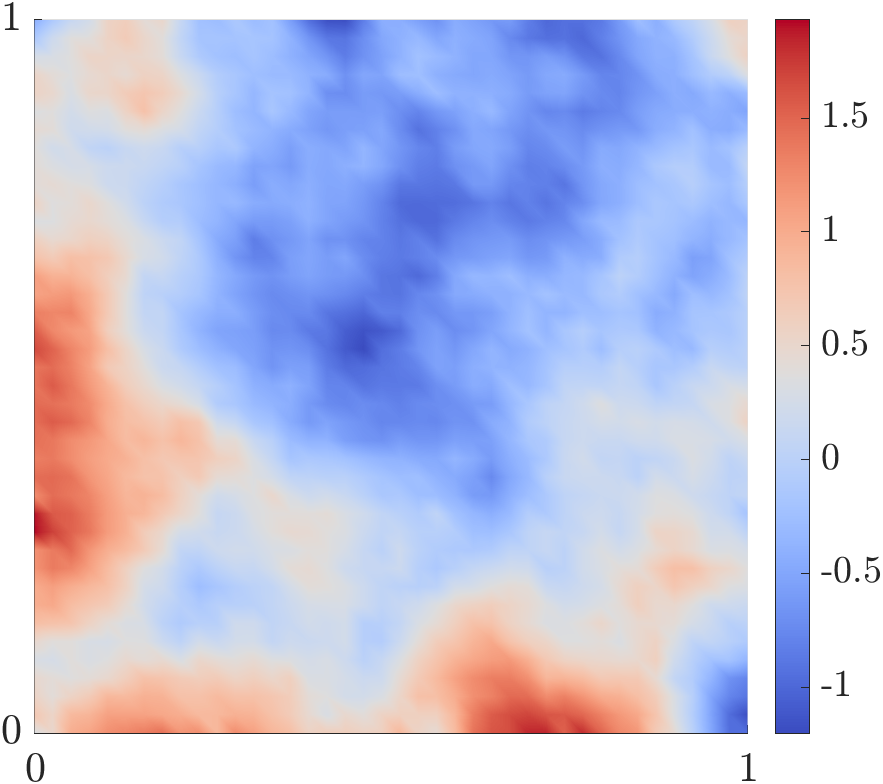}
\hfill
\includegraphics[height=0.27\textwidth]{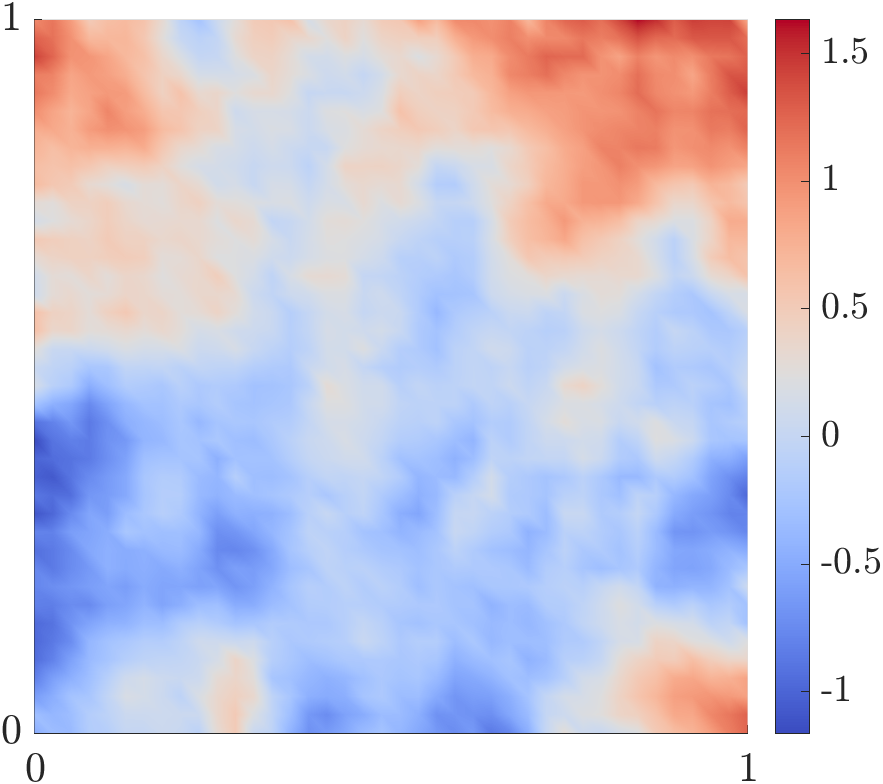}\\
$\parb$ \hrulefill\vspace{5pt}

\includegraphics[height=0.27\textwidth]{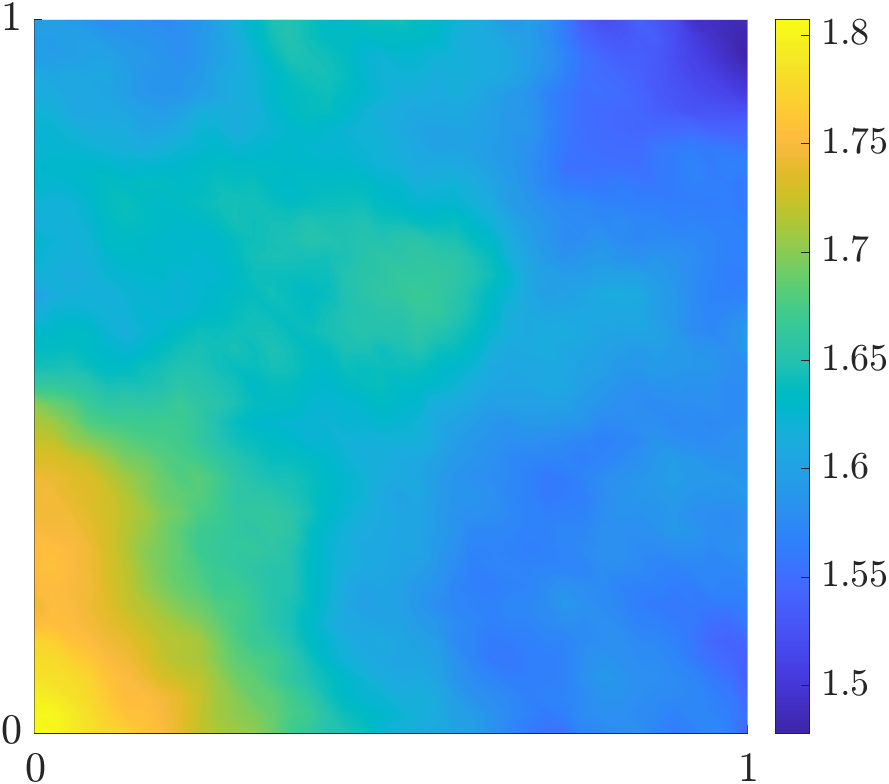}   
\hfill
\includegraphics[height=0.27\textwidth]{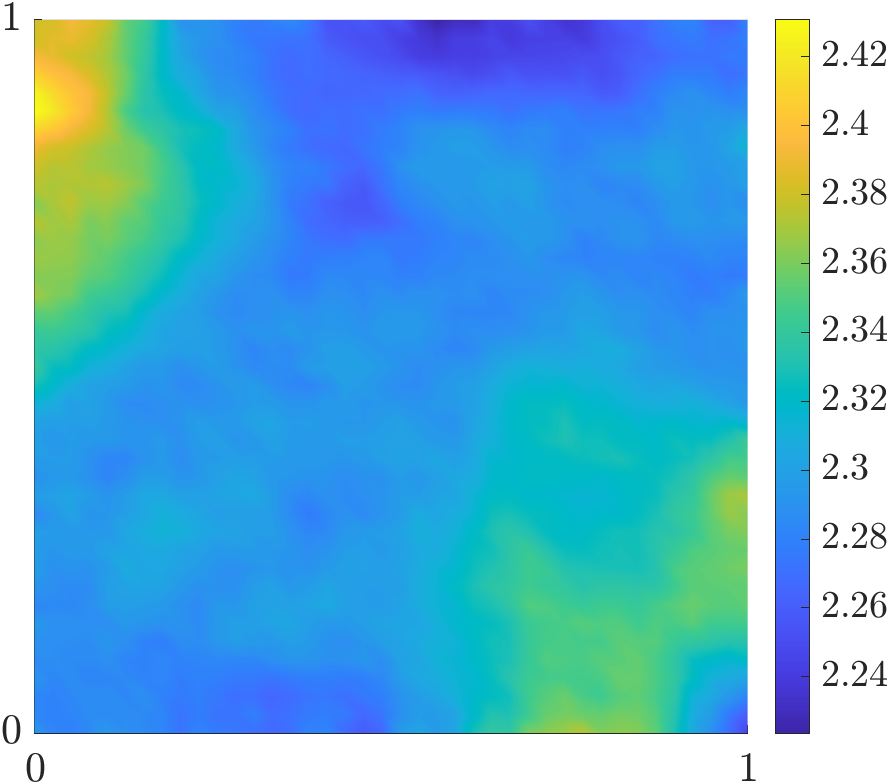}
\hfill
\includegraphics[height=0.27\textwidth]{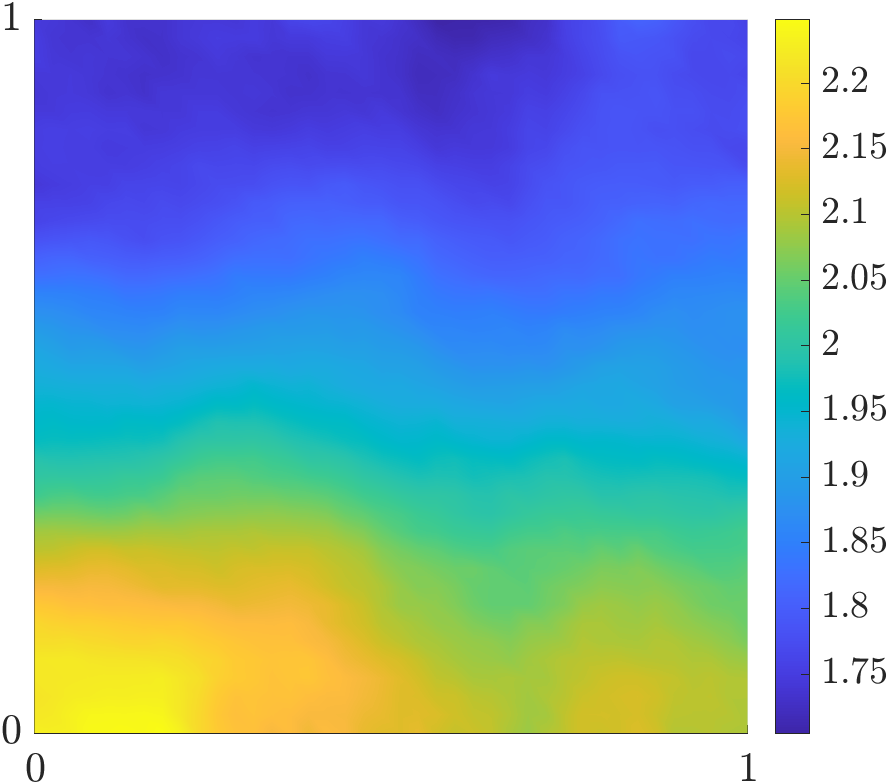}
    \caption{Example 2: Three prior samples of $\parm$ (top) and $\parb$ (bottom).}
    \label{fig: priorSampsEx1}
\end{figure}

We set $\Omega=(0,1)\times(0,1)$ and solve both the accurate forward model and the surrogate model using a Galerkin finite element method
with piecewise linear continuous Lagrange basis functions. Specifically, we consider a regular triangular mesh of the domain consisting of 1600 nodes and 3042 triangular elements which coincides with the mesh used to discretise both the parameters $\parm$ and $\parb$ (thus $n=q=1600$). The prior mean for  $\parm$ and $\parb$ are set to $\parm_0=0$ and $\parb_0=2$, respectively, while the covariance operators (see Equation (\ref{eq: priorcovs})) are set by taking $\gamma_\parm=200$, $\gamma_\parb=100$, $\kappa_\parm=4$, $\kappa_\parb=100$, and $\bs{\Theta}_\parm=\bs{\Theta}_\parb=\bs{I}$. Three prior samples of $\parm$ and $\parb$  are shown in Figure \ref{fig: priorSampsEx1}.

\begin{figure}[t!]
\includegraphics[height=0.27\textwidth]{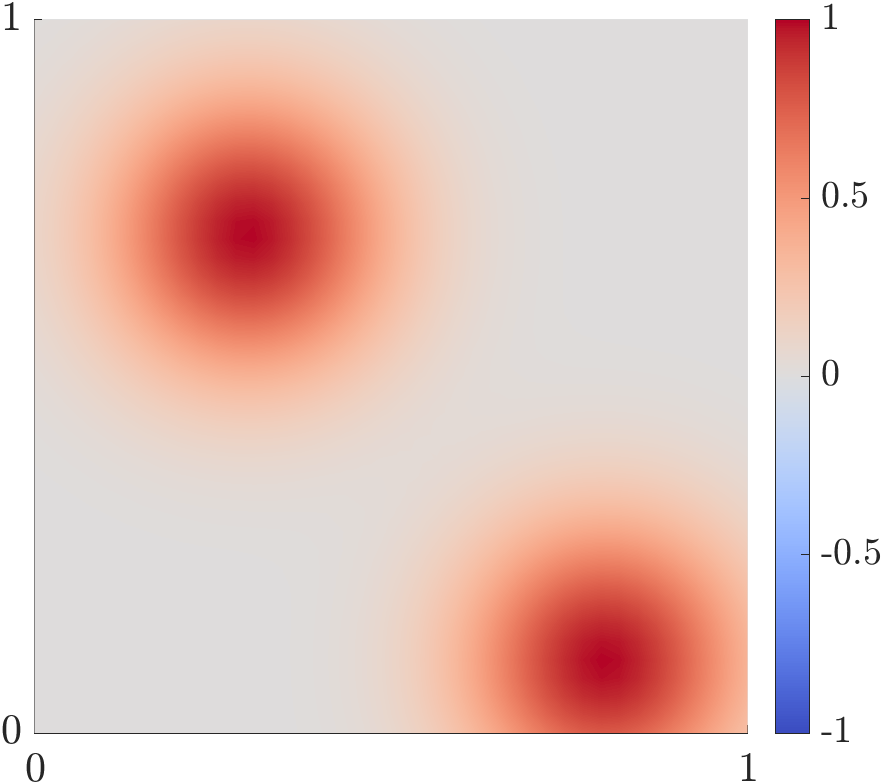}   
\hfill
\includegraphics[height=0.27\textwidth]{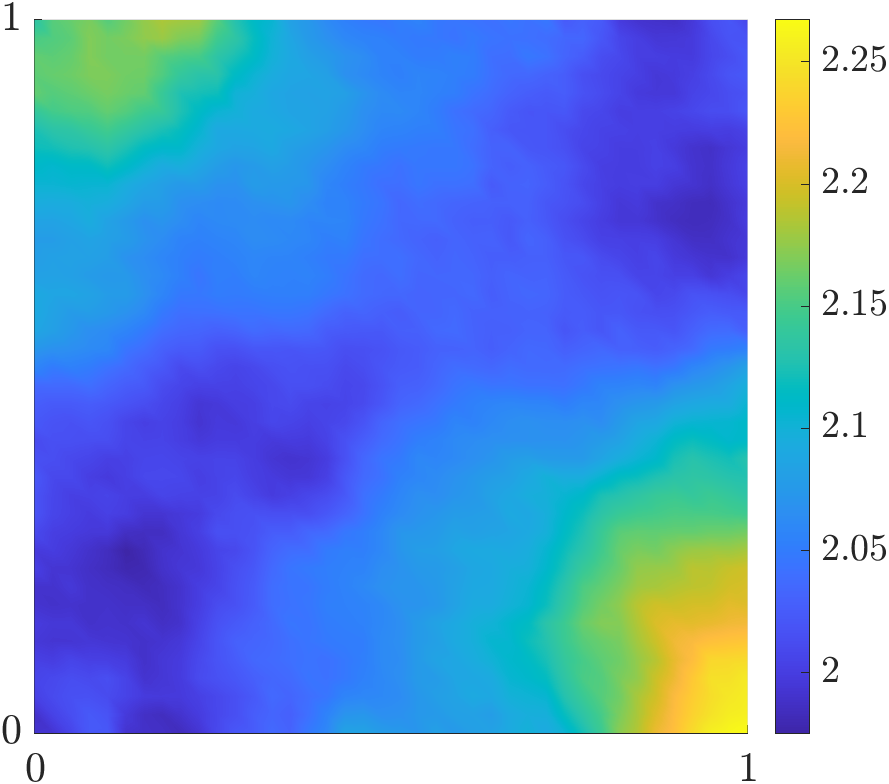}
\hfill
\includegraphics[height=0.27\textwidth]{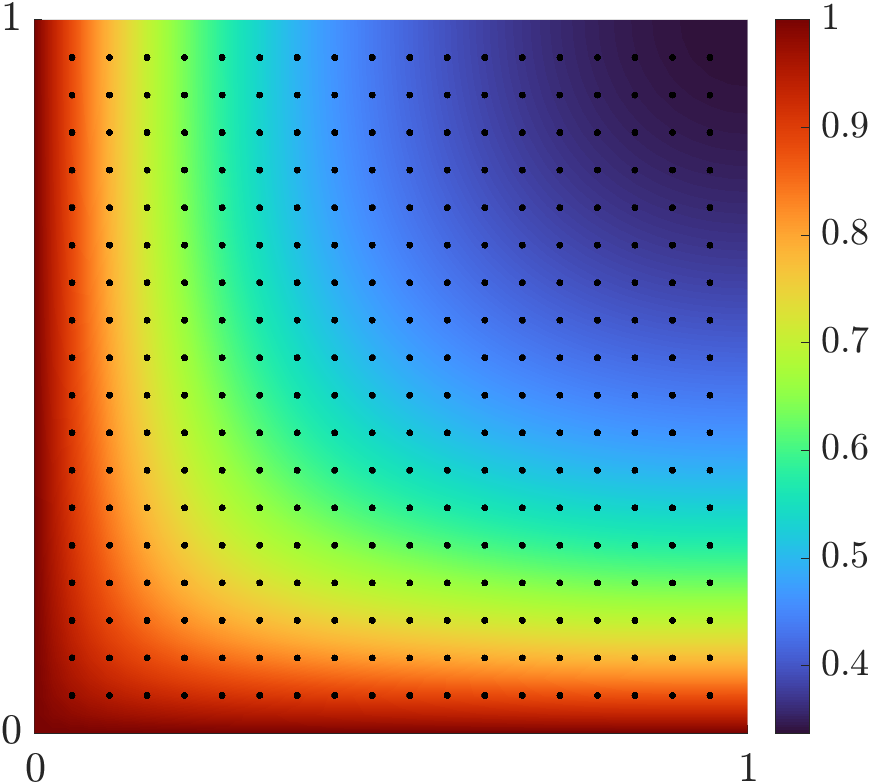}
    \caption{Example 2: The true parameters $\parm_{\rm true}$ (left) and $\parb_{\rm true}$ (centre) and the resulting solution to (\ref{eq: PAT1}) $u_{\rm true}$ (right) with measurement locations shown as black dots.}
    \label{fig: truthsEx1}
\end{figure}

The Dirichlet boundary condition for the (accurate and surrogate) forward models is set as $g=1$, with $\Gamma_{\rm D}:=(\{0\}\times [0,1])\cup( [0,1]\times\{0\})$ and $\Gamma_{\rm N}:=(\{1\}\times [0,1])\cup( [0,1]\times\{1\})$. To solve the nonlinear accurate forward model, a Newton method is used with an initial guess for $u$ computed as the solution of the linear forward model resulting from ignoring the nonlinear term, i.e., the initial guess for $u$ is the solution to the surrogate forward model. Convergence is established when the norm of the gradient reduces by a factor of $10^{10}$, typically requiring between about 4 and 6 iterations. Data consists of noisy pointwise measurements of $u$ at 324 points in a regular (18 by 18) grid throughout the domain, see Figure \ref{fig: truthsEx1}. The noise corrupting the measurements is of the form $\e\sim\mathcal{N}(\bs{0},\delta_e^2\bs{I})$, where $\delta_e=(\max(\mathcal{B}u)-\min(\mathcal{B}u))\times 1/100$, that is, the noise level is 1\% of
the range of the noiseless measurements. The true parameter of interest $\parm_{\rm true}$, true auxiliary parameter $\parb_{\rm true}$ as well as the resulting solution $u$ from which the (noiseless) data is extracted to carry out inversions are shown in Figure \ref{fig: truthsEx1}. Discussion of the BAE convergence as well as the inversion results for this example are presented in the next section.

\section{Results}\label{sec: AllRes}
To evaluate the performance of the proposed control variate approach, we first compare the convergence of the approximation error mean and covariance using standard MC estimation and the control variate-enhanced MC for both examples. We then investigate the accuracy of posterior estimates (i.e., respective MAP estimates and approximate posterior covariance matrices).

\subsection{Convergence of the approximation error statistics}\label{sec: Res1} 

To assess the accuracy and efficiency of the control variate approach, we compute reference estimates of the approximation error mean and covariance using standard MC with $N=10^5$ samples. These serve as our ``ground truth" for evaluating convergence behavior. All convergence results are averaged over 100 independent MC trials to account for variability due to random sampling.

\begin{figure}[t!]
    \centering
    \includegraphics[height=0.32\textwidth]{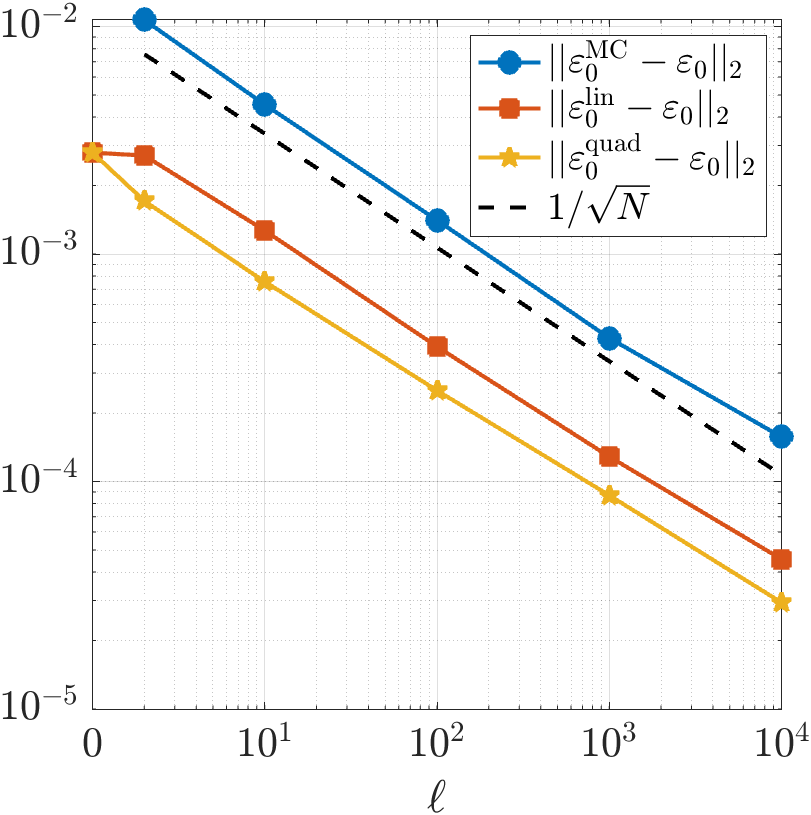}
    \hfill
            \includegraphics[height=0.32\textwidth]{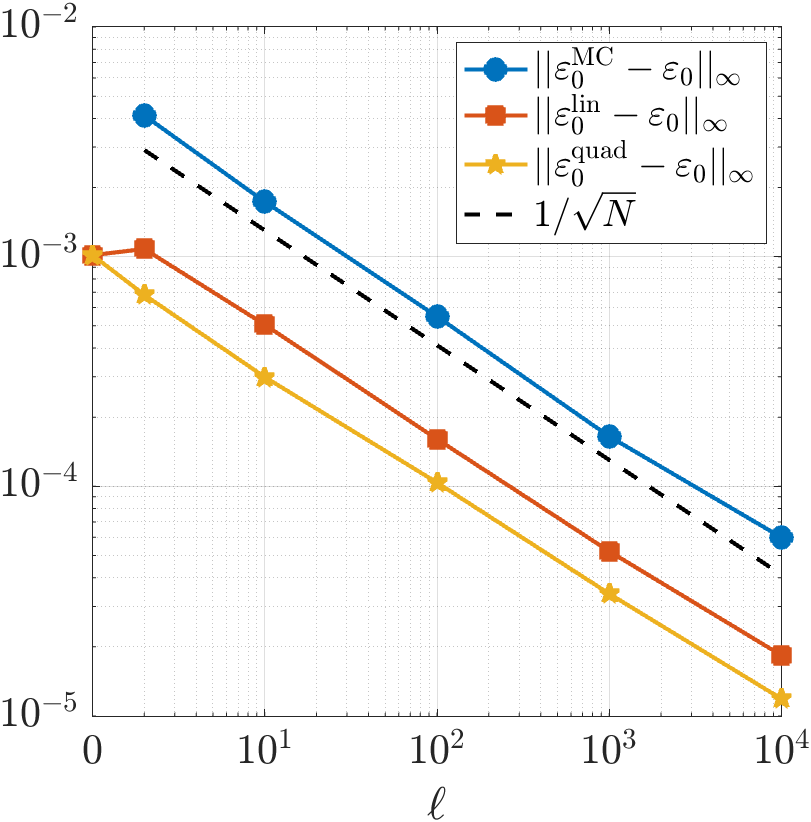}
    \hfill\includegraphics[height=0.32\textwidth]{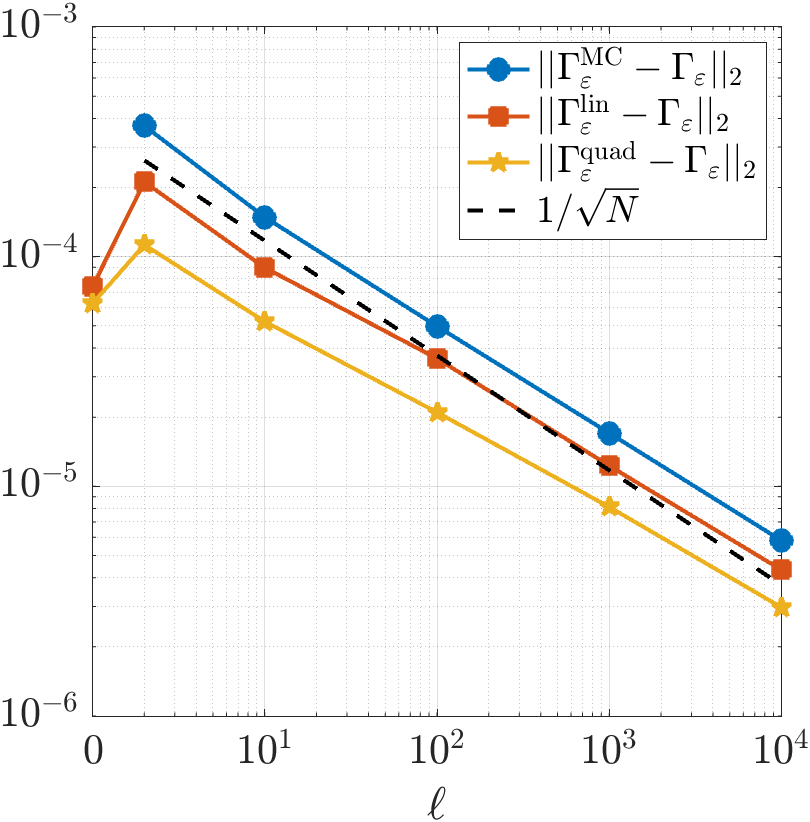}
        
    \caption{Example 1: Convergence of the mean of the approximation errors in the $\ell^2$-norm (left) and the $\ell^{\infty}$-norm (center) as well as the convergence of the covariance of the approximation errors (right). The errors for the standard MC, linear Taylor approximation control variable, and quadratic Taylor approximation control variable are shown in blue, red, and yellow respectively. Also shown is the asymptotic convergence rate $1/\sqrt{N}$ using a dashed line.}
    \label{fig: ex2Conv}
\end{figure}

Figure \ref{fig: ex2Conv} illustrates convergence of the approximation error mean in both $\ell_2$ and $\ell_\infty$ norms, and of the Frobenius norm of the error in the covariance for the first numerical example (linear diffusion problem). When using a linear Taylor approximation as a control variable, the convergence of the mean improves significantly, while use of the quadratic Taylor approximation as a control variable improves convergence further. Notably, approximately 10 samples using the linear control variable approach achieves comparable accuracy to approximately 100 samples with standard MC, while standard MC requires approximately 300 samples to achieve the same level of accuracy as 10 samples using the quadratic control variate approach.
Covariance convergence is also improved, though to a lesser extent, the quadratic control variable approach requiring approximately an order of magnitude less samples for the same accuracy as the standard MC approach.

\begin{figure}[t!]
    \centering
    \includegraphics[width=0.3\textwidth]{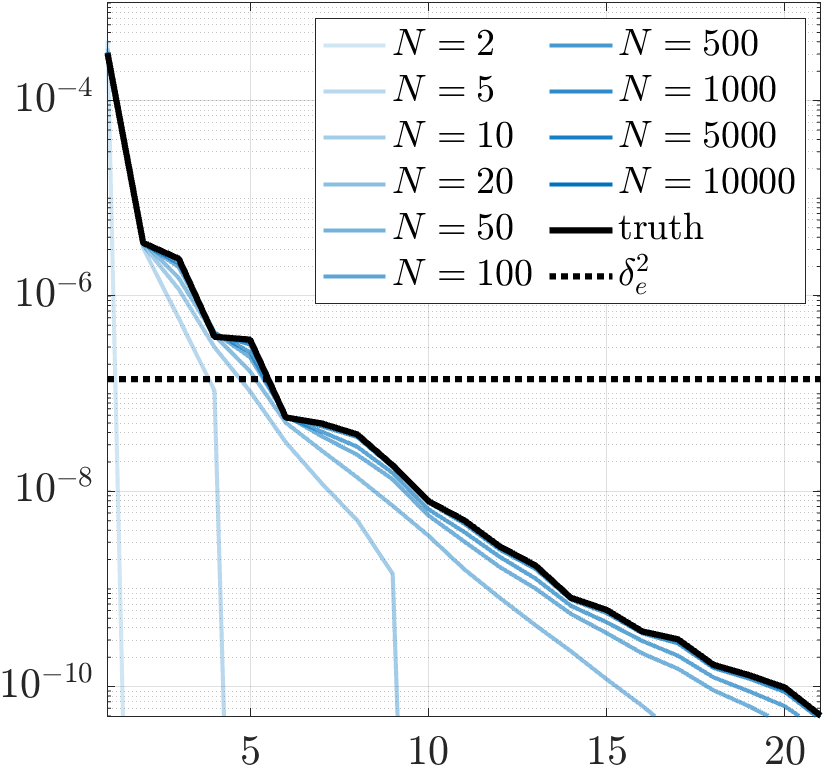}
    \hfill
        \includegraphics[width=0.3\textwidth]{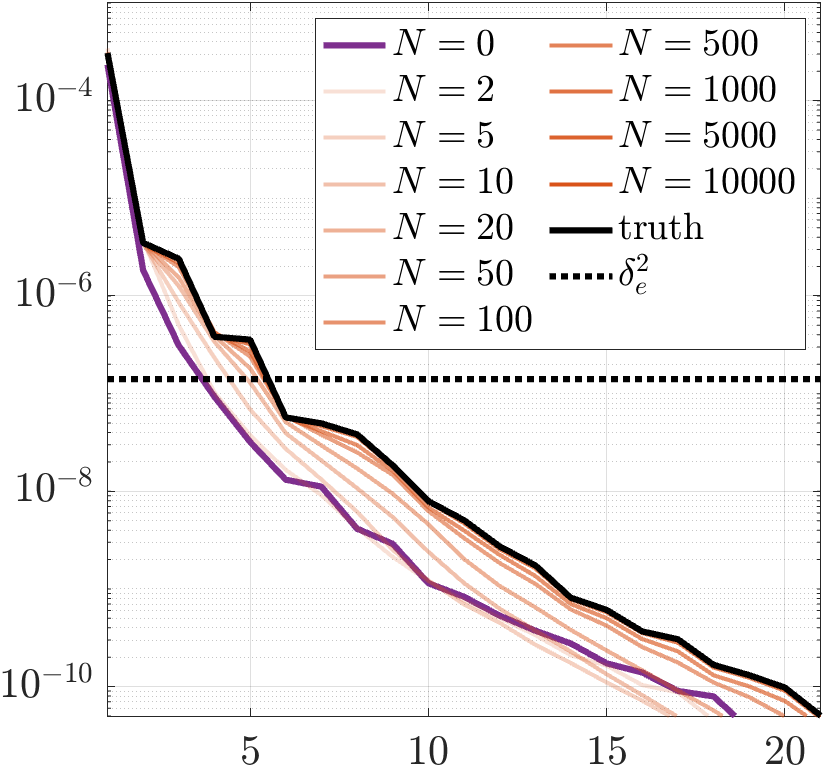}
    \hfill
        \includegraphics[width=0.3\textwidth]{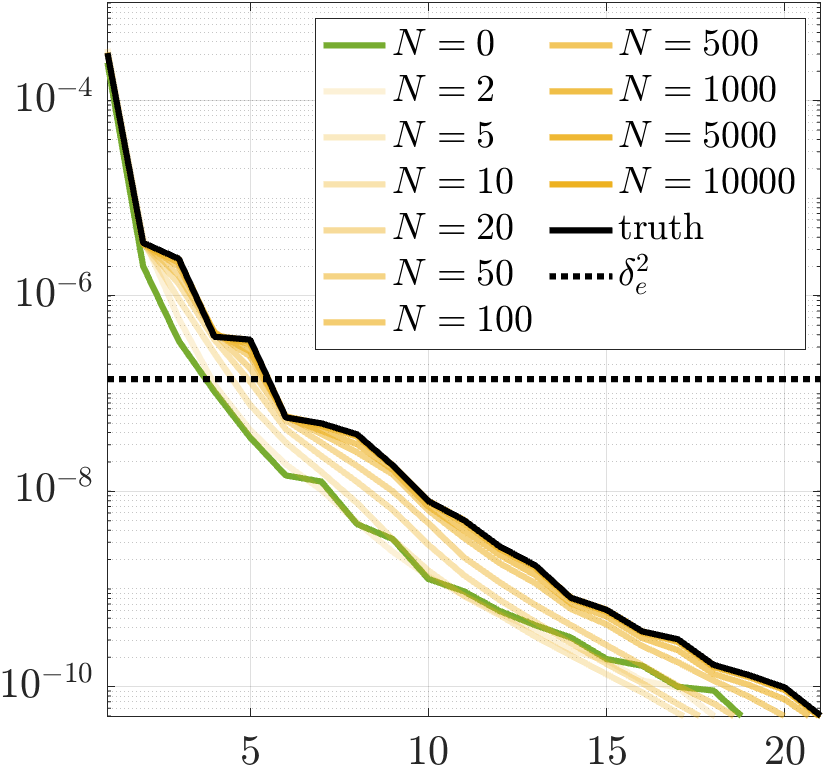}
        
    \caption{Example 1: Convergence of the spectrum of $\bs{\Gamma}_\eps$ using standard MC (left) linear CV (centre) and quadratic CV  (right) using $N\in\{0,2,5,10,20,50,100,500,1000,5000,10000\}$ samples, dotted line is the noise level variance $\delta_e^2$.}
    \label{fig: ex2Spec}
\end{figure}

To gain additional insight into the convergence of the covariance, in Figure~\ref{fig: ex2Spec} we compare the spectral decay of the estimated covariance matrix $\bs{\Gamma}_\eps$ using standard MC and when using the first- and second-order Taylor approximation as control variables. The leading eigenvalues (those above $10^{-8}\approx0.1\times\delta^2$), which are the most important to account for approximation errors (see Remark \ref{rem: noisepre}), appear to stabilise more rapidly under the control-variate schemes. Specifically, the plots indicate that to capture the behaviour of the leading eigenvalues we need $\mathcal{O}(10^1)$ approximation error samples using the control-variate approaches (the linear CV approach requiring slightly more than the quadratic CV) and  $\mathcal{O}(10^2)$ samples when using standard MC.

\paragraph{Sample-free approximation errors}
In Figures \ref{fig: ex2Conv} and \ref{fig: ex2Spec} we also show errors in the mean and covariance and the error in the spectrum of the approximation error covariance, respectively, when using the sample-free approach (see Section~\ref{sec:sampFree}). The error in the mean (in the $\ell^2$- and $\infty$-norm) found using the sample-free approximation is approximately the same as that of standard MC using $N\approx30$ samples. Moreover, in terms of the error in the covariance, the sample-free approximation attains the same accuracy as using $N\approx50$ standard MC samples, with the quadratic sample-free approximation being slightly better than the linear sample-free approximation. It is interesting to note the behaviour of the errors going from the sample-free approximation to using $N=2$ samples using the Taylor approximations as control variables. This phenomenon underlies the robustness of the sample-free approach, compared to the use of ``too few'' MC samples. 

\begin{figure}[t!]
    \centering
    \includegraphics[height=0.32\textwidth]{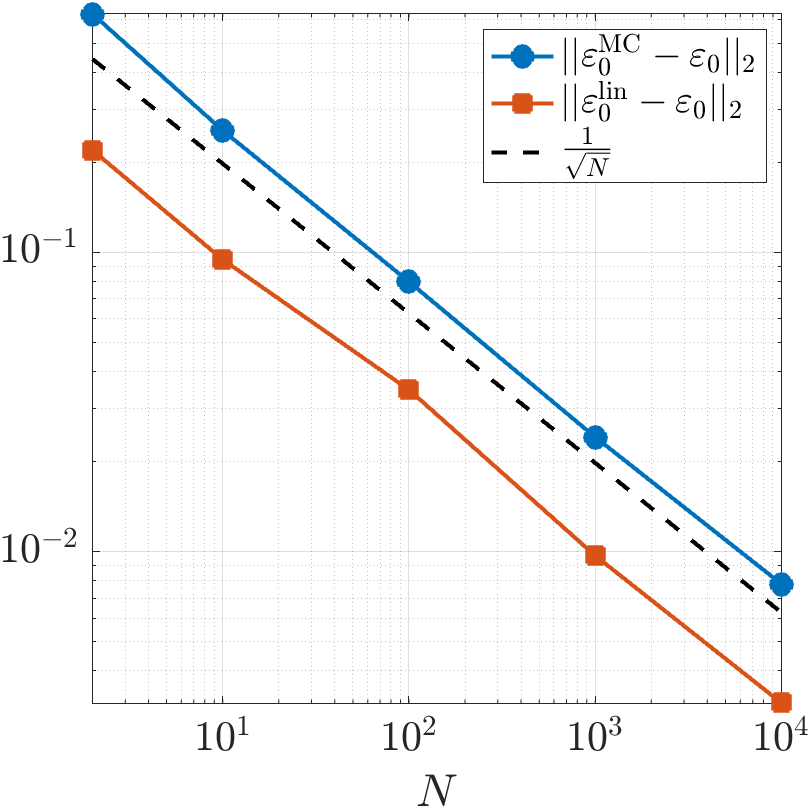}
    \hfill
     \includegraphics[height=0.32\textwidth]{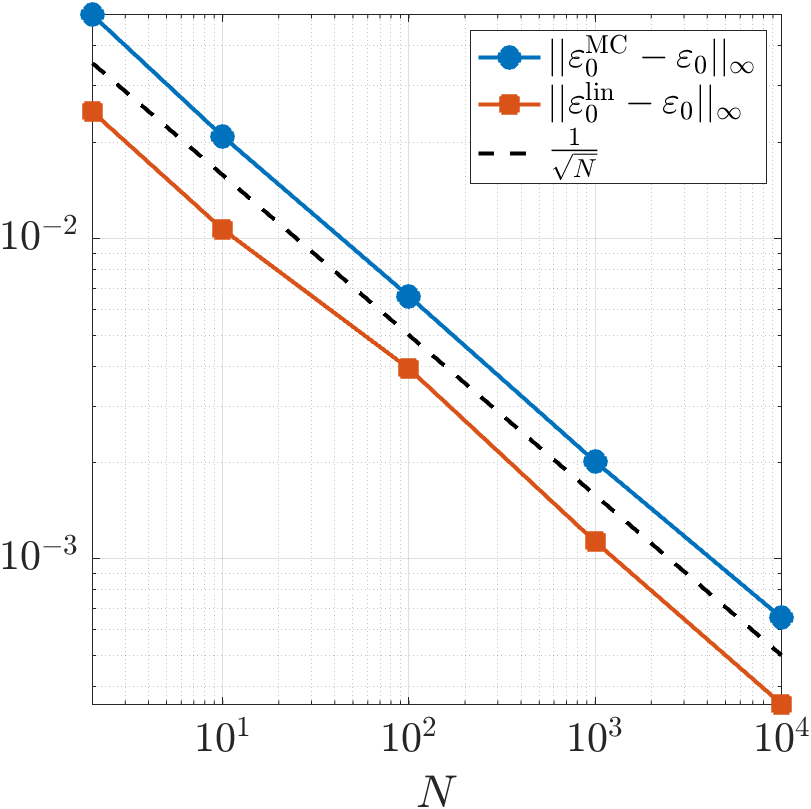}
    \hfill\includegraphics[height=0.32\textwidth]{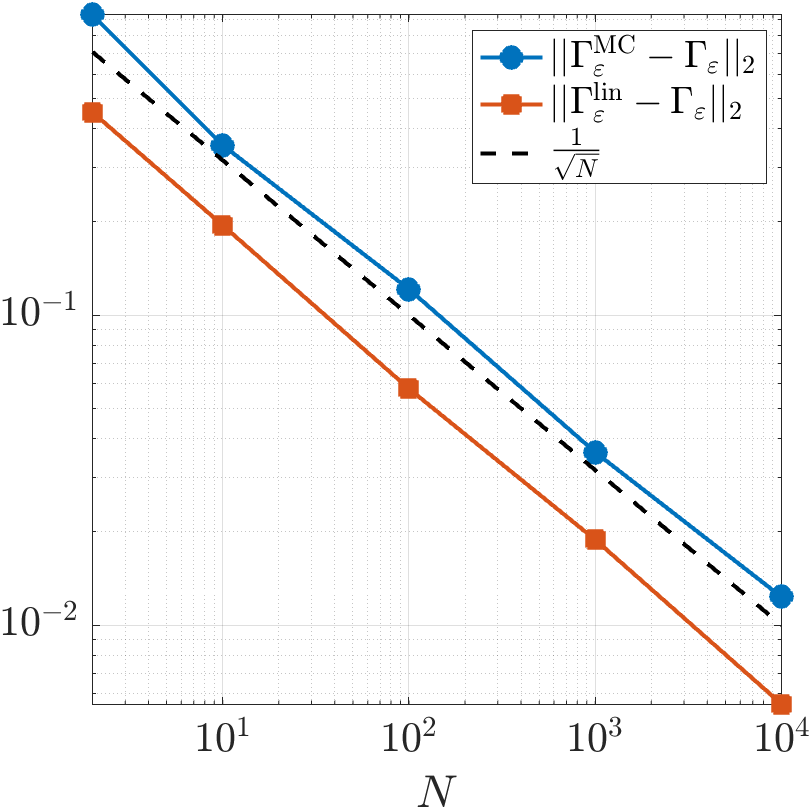}
        
    \caption{Example 2: Convergence of the mean of the approximation errors in the $\ell_2$-norm (left) and the $\ell_{\infty}$-norm (center) as well as the convergence of the covariance of the approximation errors (right).}
    \label{fig: ex1Conv}
\end{figure}

\begin{figure}[t!]
    \centering
    \includegraphics[width=0.3\textwidth]{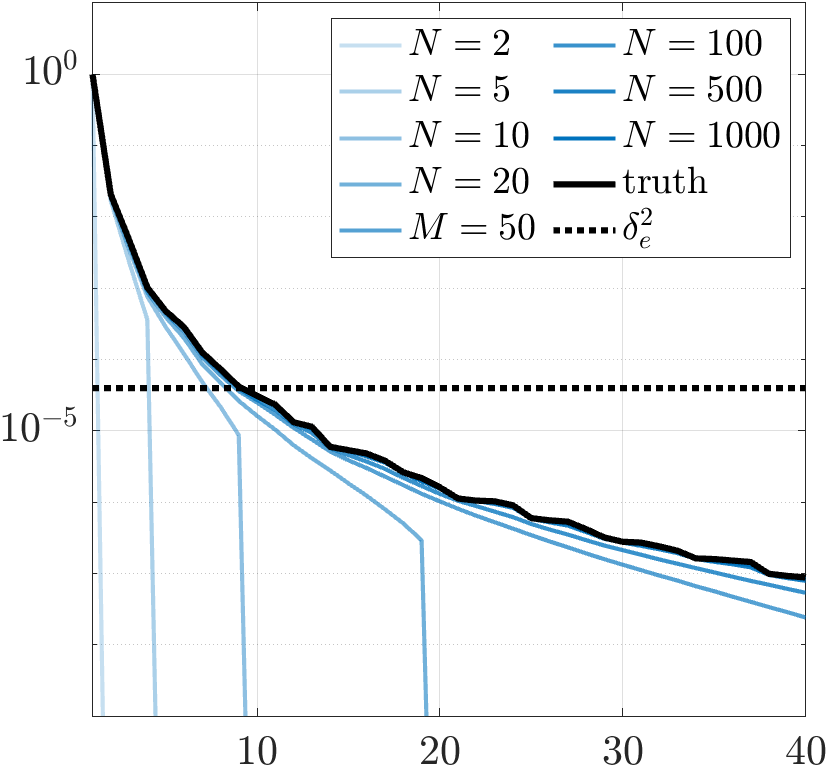}    
    \hspace{2cm}
    \includegraphics[width=0.3\textwidth]{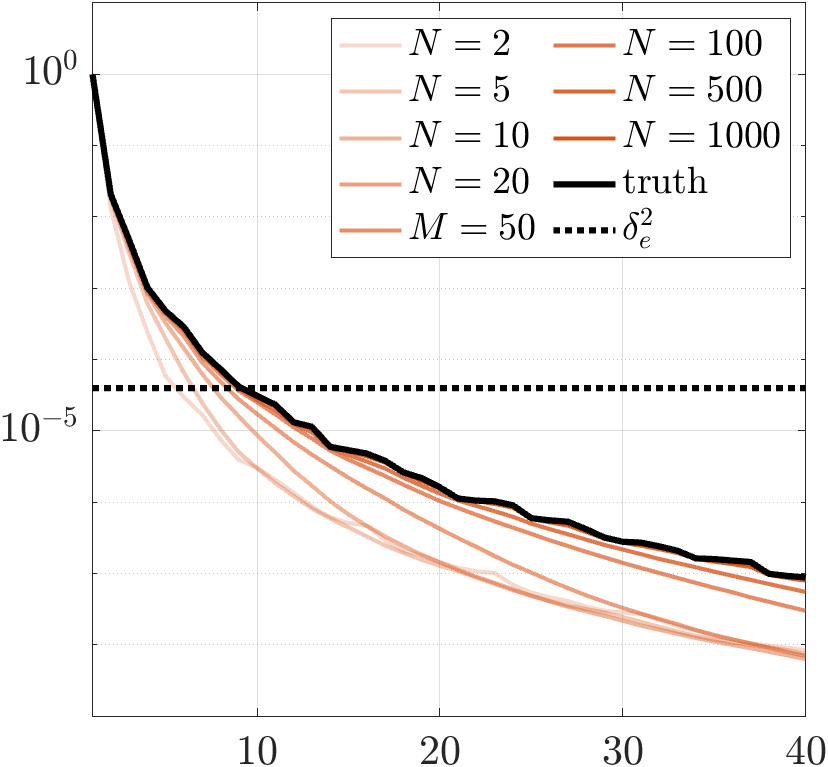}
        
    \caption{Example 2: Convergence of the spectrum of $\Gamma_\eps$ using standard MC (left) and linear CV (right) using $N\in\{2,5,10,20,50,100,500,1000\}$ samples, dotted line is the noise level variance $\delta_e^2$}
    \label{fig: ex1Spec}
\end{figure}

Similar trends are observed in the second numerical example (nonlinear diffusion problem). As shown in Figures \ref{fig: ex1Conv}, use of the linear control variable again yields marked improvement in mean estimation, as well as in the covariance, though to a lesser extent. The benefit is also visible in the stabilisation of the leading eigenvalues in Figure~\ref{fig: ex1Spec}.

\subsection{Posterior estimates}\label{sec: postests}
In this section we investigate how the use of control variate-enhanced estimation of the approximation error statistics influences the resulting (approximate) posterior distribution of the parameter of interest. Specifically, we assess the accuracy of the resulting marginal posterior 
$\pi(\bs{\parm}|\data)$ when the approximation error is estimated using either standard MC or the proposed control variate-enhanced MC methods, for both numerical examples. We restrict ourselves to comparing the respective Gaussian approximations to the marginal posterior, each of which are of the form $\mathcal{N}(\bs{\parm}_{\rm MAP}^{\rm BAE},\bs{\G}_{\parm|d}^{\rm BAE})$, with $\bs{\parm}_{\rm MAP}^{\rm BAE}$ and $\bs{\G}_{\parm|d}^{\rm BAE}$ defined in (\ref{eq:BAE_map}) and (\ref{eq:BAE_cov}), respectively. The difference between the approximations stems from how the statistics of the approximation errors are computed, e.g., using standard MC sampling or the proposed control variate-based approach. 

To assess the resulting approximate Gaussian posteriors, we compare both the MAP estimates and the squared 2-Wasserstein (Kantorovich–Rubinstein) distance \cite{villani2008optimal,villani2021topics} between each of approximate posterior distributions and a reference distribution. For both examples the reference distribution is the Gaussian
approximation to the marginal posterior using  standard MC using $N=10^5$ approximation error samples, which we denote by $\hat{\pi}^{\rm BAE}(\bs{\parm}|\data)=\mathcal{N}(\hat{\bs{\parm}}^{\rm BAE}_{\rm MAP},\hat{\bs{\G}}^{\rm BAE}_{m|d})$. The squared 2-Wasserstein distance for each of the approximate posterior distributions from 
$\hat{\pi}(\bs{\parm}|\data)$ as 
\begin{align}
   \label{eq:Wasser}W_2^2\left(\pi^\chi(\bs{\parm}|\data)\right)&=\norm[\bs{M}_\parm]{\bs{\parm}^\chi_{\rm MAP}-\hat{\bs{\parm}}_{\rm MAP}}^2+\tr\left(\bs{\G}^\chi_{m|d}+\hat{\bs{\G}}_{m|d}-2\left(\hat{\bs{\G}}_{m|d}^{\frac{1}{2}}\bs{\G}^\chi_{m|d} \hat{\bs{\G}}_{m|d}^{\frac{1}{2}}\right)^{\frac{1}{2}}\right), 
\end{align}
where $\chi\in\{\text{MC},\text{lin},\text{quad}\}$ denotes which method (Monte Carlo (MC), linear Taylor series as CV (lin), or quadratic Taylor series as CV (quad)) is used to compute the statistics of the approximation errors. To account for the stochastic nature of the problem we compute the expected value (with respect to the data) of the error in the MAP estimates and squared 2-Wasserstein distance. Specifically, we first compute the approximation errors using MC and the control variate-enhanced MC using 50 different random seeds. We then generate 20 different sets of data, i.e., we take 20 samples of $(\bs\parm^{(i)},\bs\parb^{(i)})\sim\pi(\bs\parm,\bs\parb)$ and 20 samples of $\bs e^{(i)}\sim\pi(\bs e)$ and compute $\data^{(i)}=\afwd(\bs\parm^{(i)},\bs\parb^{(i)})+\bs e^{(i)}$, with $i=1,2,\dots,20$. From each of these 20 sets of data we compute the Laplace approximation to the posterior using the BAE approach with each of the MC approaches (for each of the 50 different seeds) and compare these to the MAP estimate and the Laplace approximation of the reference distribution $\hat{\pi}^{\rm BAE}(\bs{\parm}|\data)$.

\begin{figure}[t!]
    \centering
\includegraphics[width=0.27\textwidth]{PRIOR1CV.png}    
\hfill
\includegraphics[width=0.27\textwidth]{BAE1CV.png}
\hfill
\includegraphics[width=0.27\textwidth]{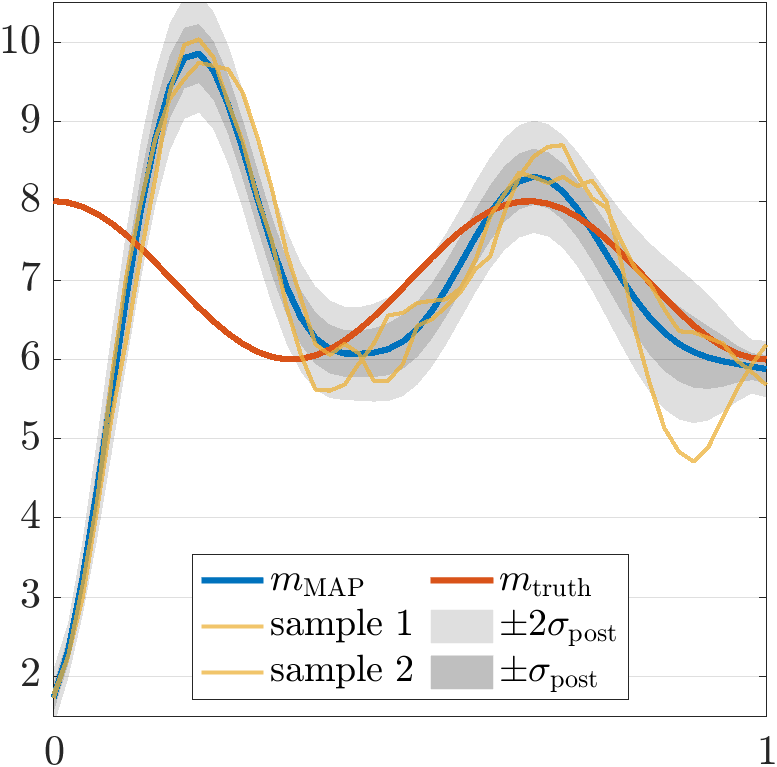}
    \caption{Posterior estimates for Example 1. The prior (left), the MAP estimate found using the reference (with $N=10^5$ samples) BAE approach   $\tilde{\parm}_{\rm MAP}^{\rm BAE} $ (centre), and the MAP estimate found ignoring the approximation errors (right). In each of the plots, the blue line indicates the mean of the respective distributions, the red line is the truth, and two samples from the respective distributions are shown in yellow. Also shown are the plus/minus one and two standard deviation intervals,  $\pm\sigma$ and $\pm2\sigma$.}
    \label{fig: MAPPOIS}
\end{figure}

\begin{figure}[!ht]
    \centering
\includegraphics[height=0.27\textwidth]{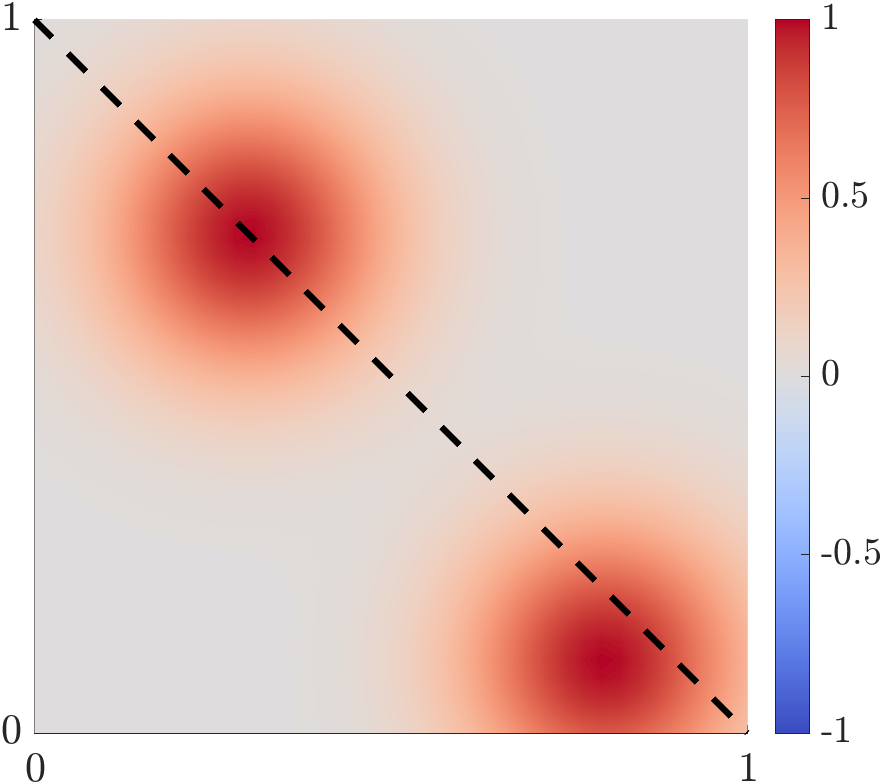}
\hfill
\includegraphics[height=0.27\textwidth]{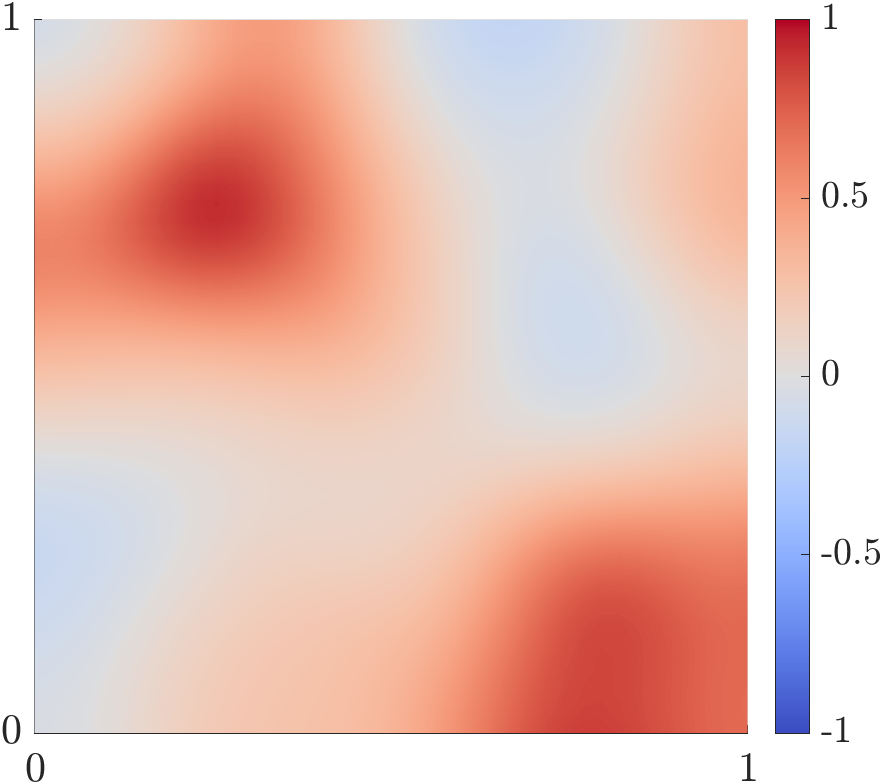}
\hfill
\includegraphics[height=0.27\textwidth]{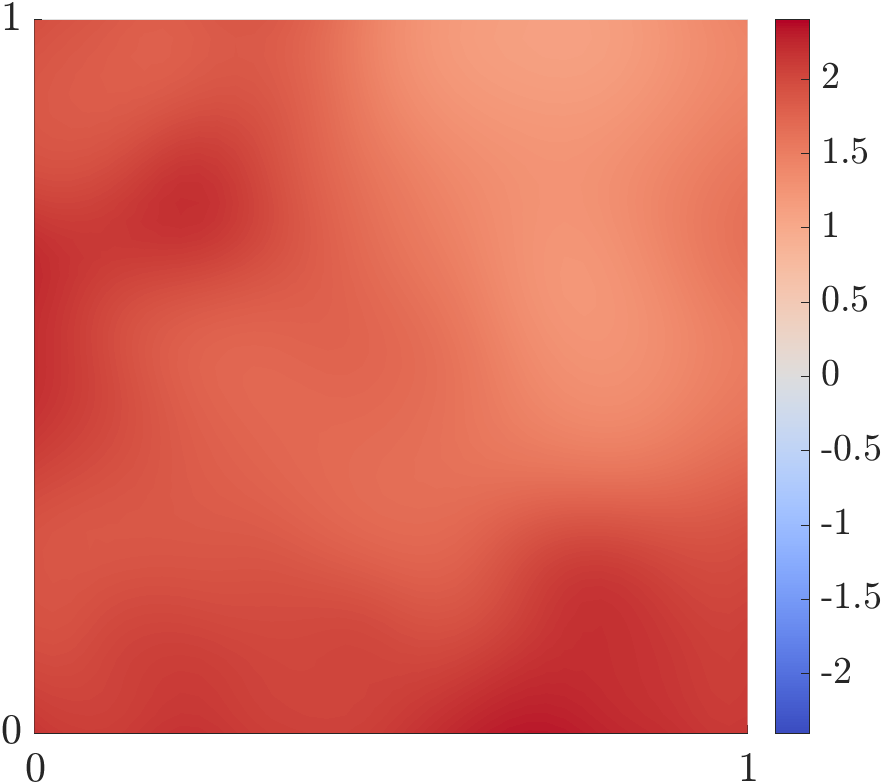}\\
\vspace{.5cm}
\includegraphics[width=0.27\textwidth]{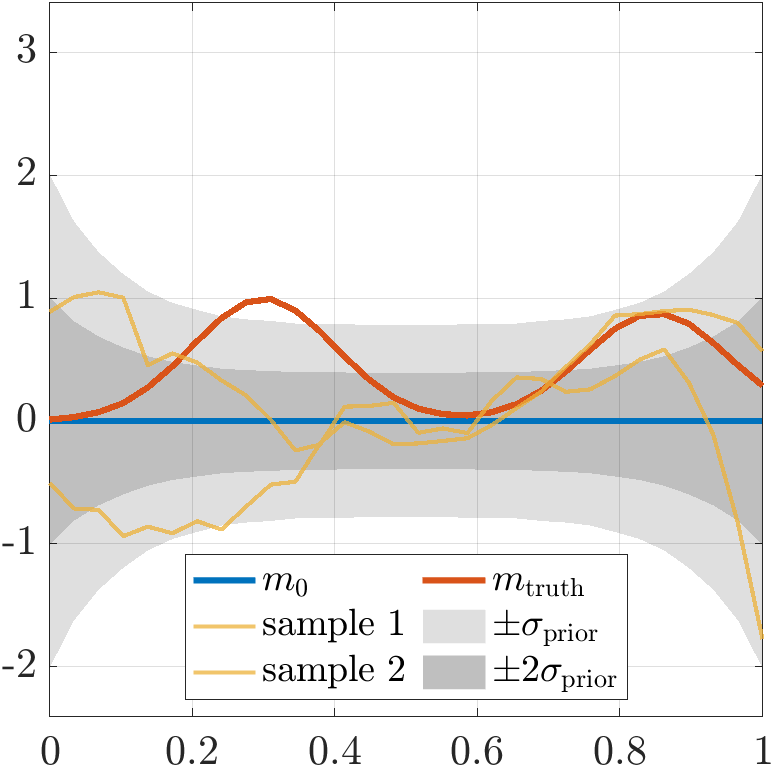}   
\hfill
\includegraphics[width=0.27\textwidth]{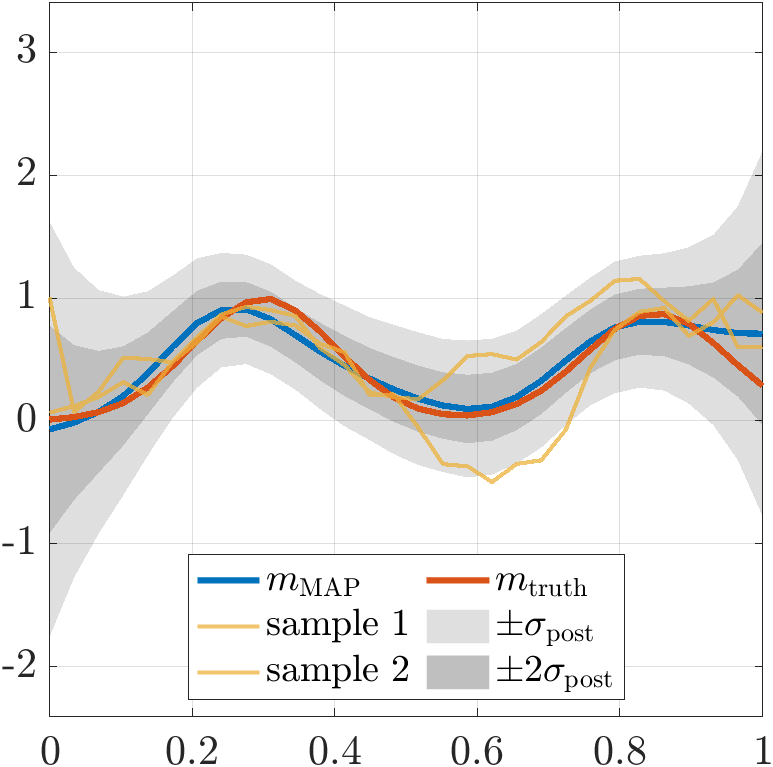}
\hfill
\includegraphics[width=0.27\textwidth]{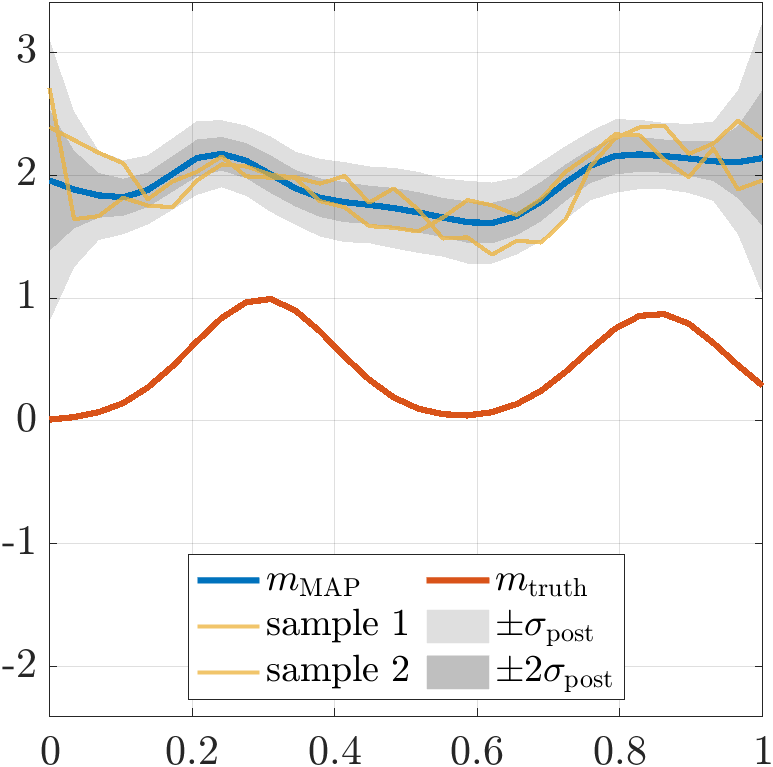}
    \caption{Posterior estimates for Example 2.  The true parameter of interest $\parm_{\rm true}$ (top left), the MAP estimate found using the reference BAE (standard MC with $N=10^5$ samples used to compute $\bs\eps_0$ and $\bs\Gamma_\eps$) approach    $\tilde{\parm}_{\rm MAP}^{\rm BAE} $ (top centre), and the MAP estimate found ignoring the approximation errors (top right), note the different colour bar. In the plot of the true parameter we show the line from $(0,1)$ to $(1,0)$ which is used to plot the one-dimensional marginal posterior plots in bottom row, where we show the prior (bottom left) and the Gaussian approximation to the posterior using the BAE approach (bottom centre) and when ignoring the approximation errors (bottom right). In each of the marginal posterior plots, the blue line indicates the mean of the respective distributions, the red line is the truth, and two samples from the respective distributions are shown in yellow. Also shown are the plus/minus one and two standard deviation intervals,  $\pm\sigma$ and $\pm2\sigma$.}
    \label{fig: MAPQPAT}
\end{figure}

The reference results, i.e., $\hat{\pi}(\bs{\parm}|\data)$, for Examples 1 and 2 are shown in Figures \ref{fig: MAPPOIS} and \ref{fig: MAPQPAT}. In these figures we also show the priors, as well as the results found when the approximation errors are ignored. In both cases, it is clear that ignoring the approximation errors leads to biased estimates and overconfident uncertainty estimates, highlighting the importance of accounting for the approximation errors in the inversion process. For comparison, in Figure~\ref{fig: REFex1} (see Appendix~\ref{sec: RefRes}) we show the posterior estimates for $\parm$ found using the true auxiliary parameter values $\parb$ in the inversions, i.e., using $\cfwd(\parm):=\afwd(\parm,\parb_{\rm true})$.

\begin{figure}[!ht]
\centering
\includegraphics[height=0.3\textwidth]{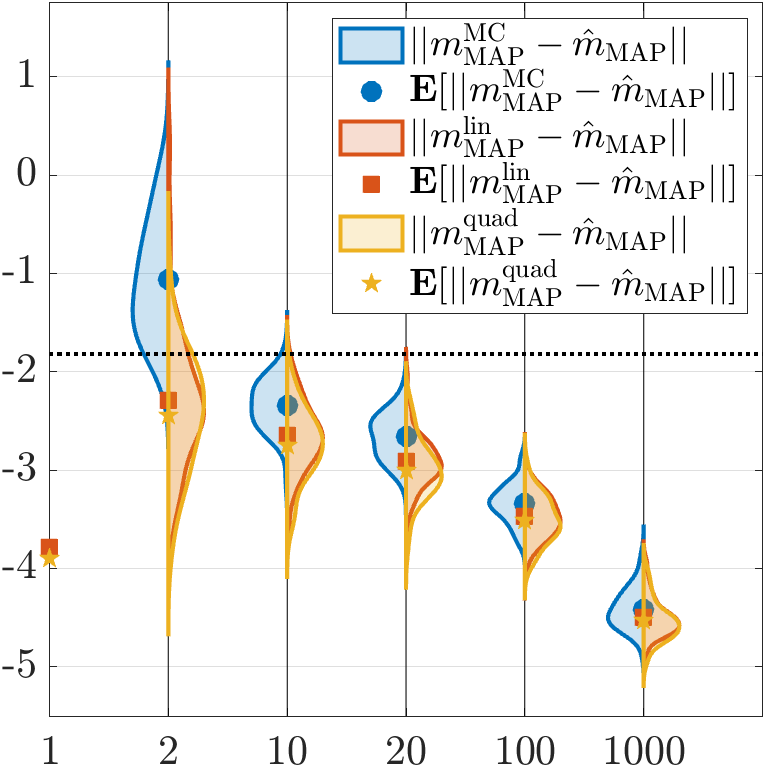} 
\hspace{2cm}
\includegraphics[height=0.3\textwidth]{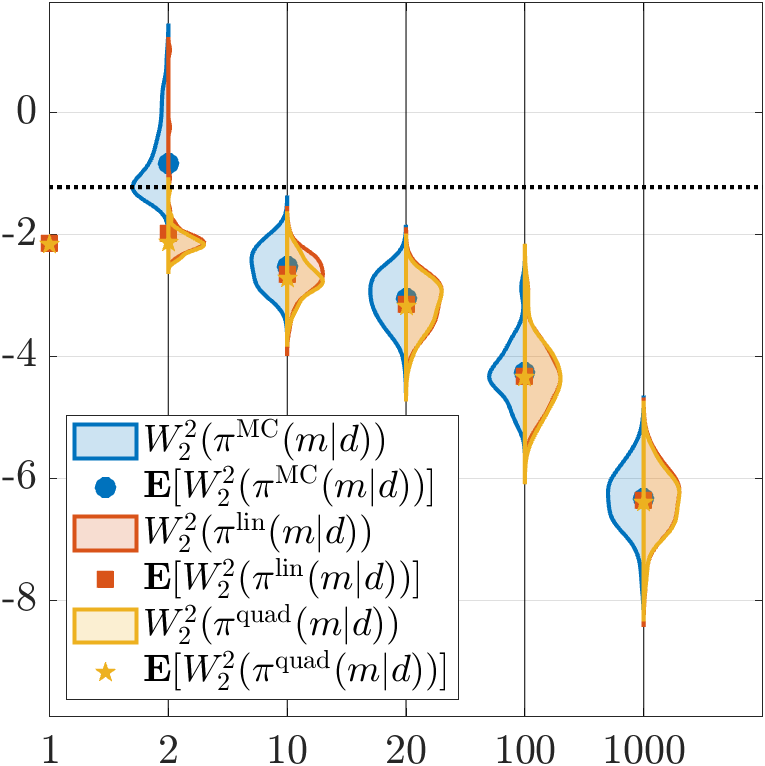}
    \caption{Example 1: Convergence of expectation (with respect to data) error in MAP estimate (left) and the squared 2-Wasserstein distance to the reference posterior (see (\ref{eq:Wasser})) (right) as a function of the number of BAE samples $N$. The results are averaged over 20 different data realisations, while 50 different random seeds are used (giving the distribution) for the BAE sampling.
    The horizontal dashed line shows the error between the reference results and the case of no approximation error, i.e., with $\parb=\parb_{\rm true}$, $\bs{\eps}_0=\bs{0}$ and $\bs{\Gamma}_\eps=\bs{0}$.}
    \label{fig: ConvPois}
\end{figure}

\begin{figure}[!ht]
\centering
\includegraphics[height=0.3\textwidth]{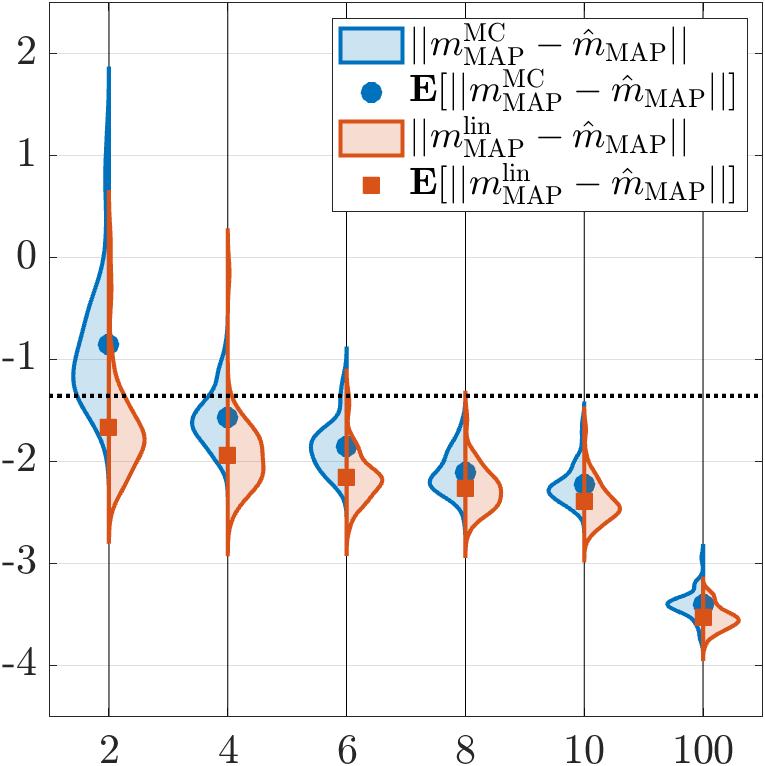}
\hspace{2cm}
\includegraphics[height=0.3\textwidth]{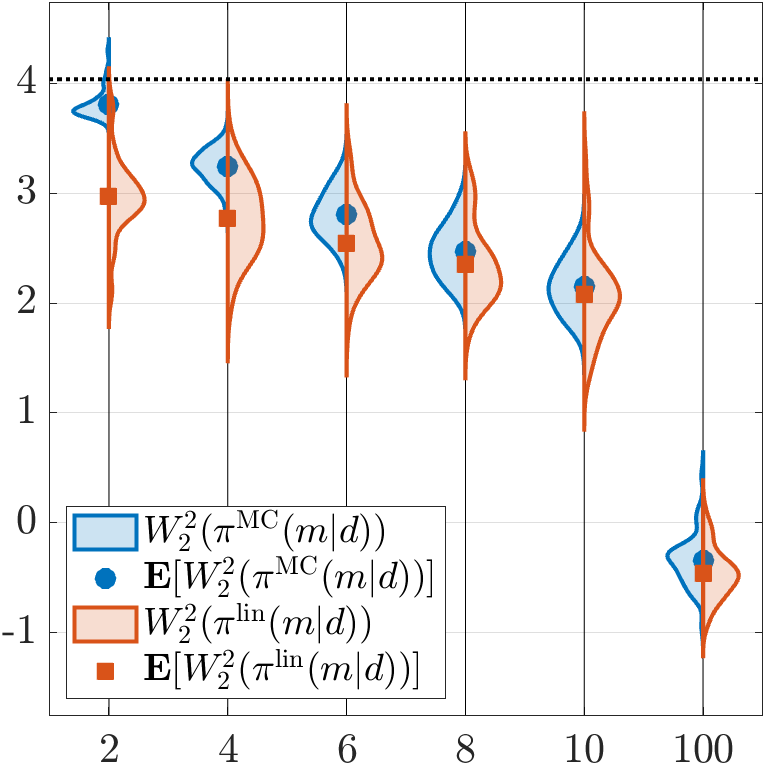}   
        \caption{Example 2: Convergence of error in the MAP estimates (left) and the squared 2-Wasserstein distance (see (\ref{eq:Wasser})) to the reference posterior (right) as a function of the number of BAE samples $N$.  The results are averaged over 20 different data realisations, while 50 different random seeds are used (giving the distribution) for the BAE sampling.The horizontal dashed line shows the error between the reference results and the case of no approximation error, i.e., with $\parb=\parb_{\rm true}$, $\bs{\eps}_0=\bs{0}$ and $\bs{\Gamma}_\eps=\bs{0}$.}
    \label{fig: ConvQPAT}
\end{figure}

In Figures~\ref{fig: ConvPois} and \ref{fig: ConvQPAT} we show the convergence of the error in the MAP estimates as well as the convergence in the squared 2-Wasserstein distance to the reference posterior for increasing sample sizes for Example 1 and Example 2, respectively. For comparison, also shown is the error in MAP and the squared 2-Wasserstein error of $\hat{\pi}^{\rm BAE}(\bs{\parm}|\data)=\mathcal{N}(\hat{\bs{\parm}}^{\rm BAE}_{\rm MAP},\hat{\bs{\G}}^{\rm BAE}_{m|d})$ to the approximate Gaussian posteriors for $\bs\parm$ when there is no approximation, i.e., when the inverse problem is solved using the ground truth value for the nuisance parameter, a comparison often considered in the BAE literature. For both examples, use of the Taylor control variates approach gives faster decay in both metrics compared to standard MC, with the quadratic approximation as control variable (used in Example 1) offering slightly better results  than the linear approximation used as the control variable. As alluded to earlier, to account for the stochastic nature of the problems, the convergence results are computed using a ``double MC'' loop using 20 different realisations of the data and 50 different random seeds for computing the BAE statistics.

\paragraph{Sample-free approximation}
Also shown in Figure~\ref{fig: ConvPois} is the error in the MAP and squared 2-Wasserstein of the sample-free approach. Similarly to the convergence of the BAE statistics reported in  Section~\ref{sec: Res1} (specifically Figure~\ref{fig: ex2Conv}), the sample-free approach can outperform the sampling-based approaches for low sample numbers $N$. Most notably, the sample-free approximation gives smaller errors in both measures than any realisation of the standard MC approach using $N=2$ samples considered here. It is also noticeable that the error in the MAP estimate for the sample-free approach is smaller (on average) than any of the sample-based approaches for up to $N=100$ samples while have $N\geq10$ samples leads to a smaller squared 2-Wasserstein. This indicates that the (approximate) posterior covariance matrix found using the sample-based approaches is more accurate than the sample-free approach (though the error is still below the error between the reference approximate posterior and the approximate posterior found when there is no approximation error). For reference, the Laplace approximation to the posterior found using the sample-free approach is shown in Figure~\ref{fig: ex1SF}.

\section{Conclusions}\label{sec: Conc}
In this work, we presented a scalable approach for reducing the computational cost of employing the Bayesian approximation error (BAE) methodology in large-scale PDE-constrained inverse problems. The core idea is to use Taylor expansions of the approximation error as control variables to accelerate Monte Carlo estimation of the error statistics. For the problems considered here, the Taylor-based control variables enabled significant variance reduction in both the mean and covariance of the approximation errors, leading to faster convergence and significantly improved efficiency in estimating the modeling error statistics.

The proposed approach retains the key advantage of the BAE framework---all samples are computed offline---while substantially reducing the number of high-fidelity forward model evaluations required during the sampling stage. Furthermore, the proposed approach is agnostic to the specific construction of the surrogate model, provided differentiability with respect to the parameters of interest. 

In addition to the control variate approach, we also explored a sample-free approximation that directly uses the analytically computed means and covariances of the Taylor approximations in place of the Monte Carlo estimates of the mean and covariance of the approximation errors. This provides a fully deterministic alternative to Monte Carlo sampling and was shown to be accurate in one of our numerical experiments.

Several extensions to the approach considered in the current paper could be used to further reduce the computational costs associated with the BAE sampling. First, it may be possible to embed the BAE sampling procedure inside an optimal experimental design (OED) framework~\cite{alexanderian2021optimalR} to develop optimal (in terms of reducing the number of samples) schemes to maximise variance reduction per solve. Somewhat related, it could be beneficial to consider employing a quasi Monte Carlo (QMC) approach to the BAE sampling, as such approaches have been shown to improve convergence rates over MC sampling~\cite{graham2011quasi,bazahica2025uncertainty,kuo2016application}. Finally, multilevel/multifidelity sampling strategies~\cite{cliffe2011multilevel,aretz2025multifidelity} could be used to further reduce cost by hierarchically combining models of varying order Taylor expansions, or more generally, of varying fidelity. In this setting, the BAE correction could be estimated at the lowest possible fidelity that preserves accuracy, with higher-fidelity models used more selectively. 

\section*{Acknowledgments}
The authors wish to thank Oliver Maclaren, Omar Ghattas, Nick Alger and Blake Christierson for several insightful discussions. This research was sponsored by the U.S. Department of Energy Office of Science
and the U.S. Air Force Office of Scientific Research. This article has been authored by an employee of
National Technology \& Engineering Solutions of Sandia, LLC under Contract No. DE-NA0003525 with
the U.S. Department of Energy (DOE). The employee owns all right, title and interest in and to the article and is solely responsible for its contents. The United States Government retains and the publisher,
by accepting the article for publication, acknowledges that the United States Government retains a non-
exclusive, paid-up, irrevocable, world-wide license to publish or reproduce the published form of this
article or allow others to do so, for United States Government purposes. The DOE will provide public
access to these results of federally sponsored research in accordance with the DOE Public Access Plan.

\bibliographystyle{siamplain}

\bibliography{bae_var}

\begin{thebibliography}{10}

\bibitem{abdallh2012bayesian}
{\sc A.~Abdallh, G.~Crevecoeur, and L.~Dupr{\'e}}, {\em A {B}ayesian approach
  for the stochastic modeling error reduction of magnetic material
  identification of an electromagnetic device}, Measurement Science and
  Technology, 23 (2012), p.~035601.

\bibitem{alexanderian2021optimalR}
{\sc A.~Alexanderian}, {\em Optimal experimental design for
  infinite-dimensional {B}ayesian inverse problems governed by {PDE}s: {A}
  review}, Inverse Problems, 37 (2021), p.~043001.

\bibitem{alexanderian2022optimal}
{\sc A.~Alexanderian, R.~Nicholson, and N.~Petra}, {\em Optimal design of
  large-scale nonlinear {B}ayesian inverse problems under model uncertainty},
  Inverse Problems, 40 (2024), p.~095001.

\bibitem{alger2020tensor}
{\sc N.~Alger, P.~Chen, and O.~Ghattas}, {\em Tensor train construction from
  tensor actions, with application to compression of large high order
  derivative tensors}, SIAM Journal on Scientific Computing, 42 (2020),
  pp.~A3516--A3539.

\bibitem{aretz2025multifidelity}
{\sc N.~Aretz, M.~Gunzburger, M.~Morlighem, and K.~Willcox}, {\em Multifidelity
  uncertainty quantification for ice sheet simulations}, Computational
  Geosciences, 29 (2025), p.~5.

\bibitem{asher2015review}
{\sc M.~J. Asher, B.~F. Croke, A.~J. Jakeman, and L.~J. Peeters}, {\em A review
  of surrogate models and their application to groundwater modeling}, Water
  Resources Research, 51 (2015), pp.~5957--5973.

\bibitem{babaniyi2021inferring}
{\sc O.~Babaniyi, R.~Nicholson, U.~Villa, and N.~Petra}, {\em Inferring the
  basal sliding coefficient field for the {S}tokes ice sheet model under
  rheological uncertainty}, The Cryosphere, 15 (2021), pp.~1731--1750.

\bibitem{bardsley2018computational}
{\sc J.~M. Bardsley}, {\em Computational Uncertainty Quantification for Inverse
  Problems: An Introduction to Singular Integrals}, SIAM, 2018.

\bibitem{bardsley2018quantitative}
{\sc P.~Bardsley, K.~Ren, and R.~Zhang}, {\em Quantitative photoacoustic
  imaging of two-photon absorption}, Journal of biomedical optics, 23 (2018),
  pp.~016002--016002.

\bibitem{bazahica2025uncertainty}
{\sc L.~Bazahica, V.~Kaarnioja, and L.~Roininen}, {\em Uncertainty
  quantification for electrical impedance tomography using quasi-{M}onte
  {C}arlo methods}, Inverse Problems, 41 (2025), p.~065002.

\bibitem{bui2013computational}
{\sc T.~Bui-Thanh, O.~Ghattas, J.~Martin, and G.~Stadler}, {\em A computational
  framework for infinite-dimensional {B}ayesian inverse problems: {P}art {I}:
  {T}he linearized case, with application to global seismic inversion}, SIAM
  Journal on Scientific Computing, 35 (2013), pp.~A2494--A2523.

\bibitem{calvetti2015artificial}
{\sc D.~Calvetti, P.~J. Hadwin, J.~P. Kaipio, and E.~Somersalo}, {\em
  Artificial boundary conditions and domain truncation in electrical impedance
  tomography. {P}art {II}: {S}tochastic extension of the boundary map}, Inverse
  Problems \& Imaging, 9 (2015), pp.~767--789.

\bibitem{calvetti2007introduction}
{\sc D.~Calvetti and E.~Somersalo}, {\em An introduction to {B}ayesian
  scientific computing: {T}en lectures on subjective computing}, vol.~2,
  Springer Science \& Business Media, 2007.

\bibitem{calvetti2023bayesian}
\leavevmode\vrule height 2pt depth -1.6pt width 23pt, {\em Bayesian scientific
  computing}, vol.~215, Springer, 2023.

\bibitem{candiani2021approximation}
{\sc V.~Candiani, N.~Hyv{\"o}nen, J.~P. Kaipio, and V.~Kolehmainen}, {\em
  Approximation error method for imaging the human head by electrical impedance
  tomography}, Inverse Problems, 37 (2021), p.~125008.

\bibitem{castello2019modeling}
{\sc D.~A. Castello and J.~P. Kaipio}, {\em Modeling errors due to {T}imoshenko
  approximation in damage identification}, International Journal for Numerical
  Methods in Engineering, 120 (2019), pp.~1148--1162.

\bibitem{chen2019taylor}
{\sc P.~Chen, U.~Villa, and O.~Ghattas}, {\em Taylor approximation and variance
  reduction for pde-constrained optimal control under uncertainty}, Journal of
  Computational Physics, 385 (2019), pp.~163--186.

\bibitem{cliffe2011multilevel}
{\sc K.~A. Cliffe, M.~B. Giles, R.~Scheichl, and A.~L. Teckentrup}, {\em
  Multilevel {M}onte {C}arlo methods and applications to elliptic {PDE}s with
  random coefficients}, Computing and Visualization in Science, 14 (2011),
  p.~3.

\bibitem{dodwell2019multilevel}
{\sc T.~J. Dodwell, C.~Ketelsen, R.~Scheichl, and A.~L. Teckentrup}, {\em
  Multilevel {M}arkov chain {M}onte {C}arlo}, Siam Review, 61 (2019),
  pp.~509--545.

\bibitem{frangos2010surrogate}
{\sc M.~Frangos, Y.~Marzouk, K.~Willcox, and B.~van Bloemen~Waanders}, {\em
  Surrogate and reduced-order modeling: {A} comparison of approaches for
  large-scale statistical inverse problems}, Large-Scale Inverse Problems and
  Quantification of Uncertainty,  (2010), pp.~123--149.

\bibitem{graham2011quasi}
{\sc I.~G. Graham, F.~Kuo, D.~Nuyens, R.~Scheichl, and I.~Sloan}, {\em
  Quasi-{M}onte {C}arlo methods for elliptic {PDE}s with random coefficients
  and applications}, Journal of Computational Physics, 230 (2011),
  pp.~3668--3694.

\bibitem{heidenreich2015bayesian}
{\sc S.~Heidenreich, H.~Gross, and M.~Bar}, {\em Bayesian approach to the
  statistical inverse problem of scatterometry: Comparison of three surrogate
  models}, International Journal for Uncertainty Quantification, 5 (2015).

\bibitem{higdon2004combining}
{\sc D.~Higdon, M.~Kennedy, J.~C. Cavendish, J.~A. Cafeo, and R.~D. Ryne}, {\em
  Combining field data and computer simulations for calibration and
  prediction}, SIAM Journal on Scientific Computing, 26 (2004), pp.~448--466.

\bibitem{huttunen2007approximation}
{\sc J.~Huttunen and J.~Kaipio}, {\em Approximation error analysis in nonlinear
  state estimation with an application to state-space identification}, Inverse
  Problems, 23 (2007), p.~2141.

\bibitem{huttunen2018approximation}
{\sc J.~M. Huttunen, J.~P. Kaipio, and H.~Haario}, {\em Approximation error
  approach in spatiotemporally chaotic models with application to
  {K}uramoto--{S}ivashinsky equation}, Computational Statistics \& Data
  Analysis, 123 (2018), pp.~13--31.

\bibitem{inglese1997inverse}
{\sc G.~Inglese}, {\em An inverse problem in corrosion detection}, Inverse
  problems, 13 (1997), p.~977.

\bibitem{KolehmainenBAE}
{\sc J.~Kaipio and V.~Kolehmainen}, {\em {Approximate marginalization over
  modelling errors and uncertainties in inverse problems}}, in {Bayesian Theory
  and Applications}, Oxford University Press, 01 2013.

\bibitem{kaipio2006statistical}
{\sc J.~Kaipio and E.~Somersalo}, {\em Statistical and computational inverse
  problems}, vol.~160, Springer Science \& Business Media, 2006.

\bibitem{kaipio2007statistical}
\leavevmode\vrule height 2pt depth -1.6pt width 23pt, {\em Statistical inverse
  problems: discretization, model reduction and inverse crimes}, Journal of
  computational and applied mathematics, 198 (2007), pp.~493--504.

\bibitem{kaipio2013bayesian}
{\sc J.~P. Kaipio and V.~Kolehmainen}, {\em Bayesian theory and applications},
  Approximate Marginalization Over Modeling Errors and Uncertainties in Inverse
  Problems,  (2013), pp.~644--672.

\bibitem{kennedy2001bayesian}
{\sc M.~C. Kennedy and A.~O'Hagan}, {\em Bayesian calibration of computer
  models}, Journal of the Royal Statistical Society: Series B (Statistical
  Methodology), 63 (2001), pp.~425--464.

\bibitem{kolehmainen2011marginalization}
{\sc V.~Kolehmainen, T.~Tarvainen, S.~R. Arridge, and J.~P. Kaipio}, {\em
  Marginalization of uninteresting distributed parameters in inverse
  problems-application to diffuse optical tomography}, International Journal
  for Uncertainty Quantification, 1 (2011).

\bibitem{kuo2016application}
{\sc F.~Y. Kuo and D.~Nuyens}, {\em Application of quasi-{M}onte {C}arlo
  methods to elliptic {PDE}s with random diffusion coefficients: {A} survey of
  analysis and implementation}, Foundations of Computational Mathematics, 16
  (2016), pp.~1631--1696.

\bibitem{law2007simulation}
{\sc A.~M. Law, W.~D. Kelton, and W.~D. Kelton}, {\em Simulation modeling and
  analysis}, vol.~3, Mcgraw-hill New York, 2007.

\bibitem{maddison2019automated}
{\sc J.~R. Maddison, D.~N. Goldberg, and B.~D. Goddard}, {\em Automated
  calculation of higher order partial differential equation constrained
  derivative information}, SIAM Journal on Scientific Computing, 41 (2019),
  pp.~C417--C445.

\bibitem{nicholson2021joint}
{\sc R.~Nicholson and M.~Niskanen}, {\em Joint estimation of {R}obin
  coefficient and domain boundary for the {P}oisson problem}, Inverse Problems,
  38 (2021), p.~015008.

\bibitem{nicholson2018estimation}
{\sc R.~Nicholson, N.~Petra, and J.~P. Kaipio}, {\em Estimation of the {R}obin
  coefficient field in a {P}oisson problem with uncertain conductivity field},
  Inverse Problems, 34 (2018), p.~115005.

\bibitem{nicholson2023global}
{\sc R.~Nicholson, N.~Petra, U.~Villa, and J.~P. Kaipio}, {\em On global normal
  linear approximations for nonlinear {B}ayesian inverse problems}, Inverse
  Problems, 39 (2023), p.~054001.

\bibitem{nissinen2010compensation}
{\sc A.~Nissinen, V.~P. Kolehmainen, and J.~P. Kaipio}, {\em Compensation of
  modelling errors due to unknown domain boundary in electrical impedance
  tomography}, IEEE transactions on medical imaging, 30 (2010), pp.~231--242.

\bibitem{peherstorfer2018survey}
{\sc B.~Peherstorfer, K.~Willcox, and M.~Gunzburger}, {\em Survey of
  multifidelity methods in uncertainty propagation, inference, and
  optimization}, Siam Review, 60 (2018), pp.~550--591.

\bibitem{petra2014computational}
{\sc N.~Petra, J.~Martin, G.~Stadler, and O.~Ghattas}, {\em A computational
  framework for infinite-dimensional {B}ayesian inverse problems, {P}art {II}:
  Stochastic {N}ewton {MCMC} with application to ice sheet flow inverse
  problems}, SIAM Journal on Scientific Computing, 36 (2014), pp.~A1525--A1555.

\bibitem{provost1992quadratic}
{\sc S.~Provost and A.~Mathai}, {\em Quadratic Forms in Random Variables:
  {T}heory and Applications}, Statistics: textbooks and monographs, Marcel
  Dekker, 1992.

\bibitem{ren2018nonlinear}
{\sc K.~Ren and R.~Zhang}, {\em Nonlinear quantitative photoacoustic tomography
  with two-photon absorption}, SIAM Journal on Applied Mathematics, 78 (2018),
  pp.~479--503.

\bibitem{rimpilainen2019improved}
{\sc V.~Rimpil{\"a}inen, A.~Koulouri, F.~Lucka, J.~P. Kaipio, and C.~H.
  Wolters}, {\em Improved {EEG} source localization with {B}ayesian uncertainty
  modelling of unknown skull conductivity}, NeuroImage, 188 (2019),
  pp.~252--260.

\bibitem{robert1999monte}
{\sc C.~P. Robert, G.~Casella, and G.~Casella}, {\em Monte Carlo statistical
  methods}, vol.~2, Springer, 1999.

\bibitem{ross2022simulation}
{\sc S.~M. Ross}, {\em Simulation}, academic press, 2022.

\bibitem{stuart2010inverse}
{\sc A.~Stuart}, {\em Inverse problems: {A} {B}ayesian perspective}, Acta
  numerica, 19 (2010), pp.~451--559.

\bibitem{tarvainen2009approximation}
{\sc T.~Tarvainen, V.~Kolehmainen, A.~Pulkkinen, M.~Vauhkonen, M.~Schweiger,
  S.~Arridge, and J.~Kaipio}, {\em An approximation error approach for
  compensating for modelling errors between the radiative transfer equation and
  the diffusion approximation in diffuse optical tomography}, Inverse Problems,
  26 (2009), p.~015005.

\bibitem{tick2019modelling}
{\sc J.~Tick, A.~Pulkkinen, and T.~Tarvainen}, {\em Modelling of errors due to
  speed of sound variations in photoacoustic tomography using a {B}ayesian
  framework}, Biomedical physics \& engineering express, 6 (2019), p.~015003.

\bibitem{villa2021hippylib}
{\sc U.~Villa, N.~Petra, and O.~Ghattas}, {\em h{IPPY}lib: {A}n extensible
  software framework for large-scale inverse problems governed by {PDE}s:
  {P}art {I}: {D}eterministic inversion and linearized {B}ayesian inference},
  ACM Transactions on Mathematical Software (TOMS), 47 (2021), pp.~1--34.

\bibitem{villani2021topics}
{\sc C.~Villani}, {\em Topics in optimal transportation}, vol.~58, American
  Mathematical Soc., 2021.

\bibitem{villani2008optimal}
{\sc C.~Villani et~al.}, {\em Optimal transport: old and new}, vol.~338,
  Springer, 2008.

\bibitem{zhu2018bayesian}
{\sc Y.~Zhu and N.~Zabaras}, {\em Bayesian deep convolutional encoder--decoder
  networks for surrogate modeling and uncertainty quantification}, Journal of
  Computational Physics, 366 (2018), pp.~415--447.

\end{thebibliography}

\appendix

\section{Reference  results}\label{sec: RefRes}
\begin{figure}[t!]
Example 1 \noindent\rule{0.15\textwidth}{0.4pt}\vspace{5pt}\hfill Example 2 \noindent\rule{0.5\textwidth}{0.4pt}\vspace{5pt}

\includegraphics[height=0.27\textwidth]{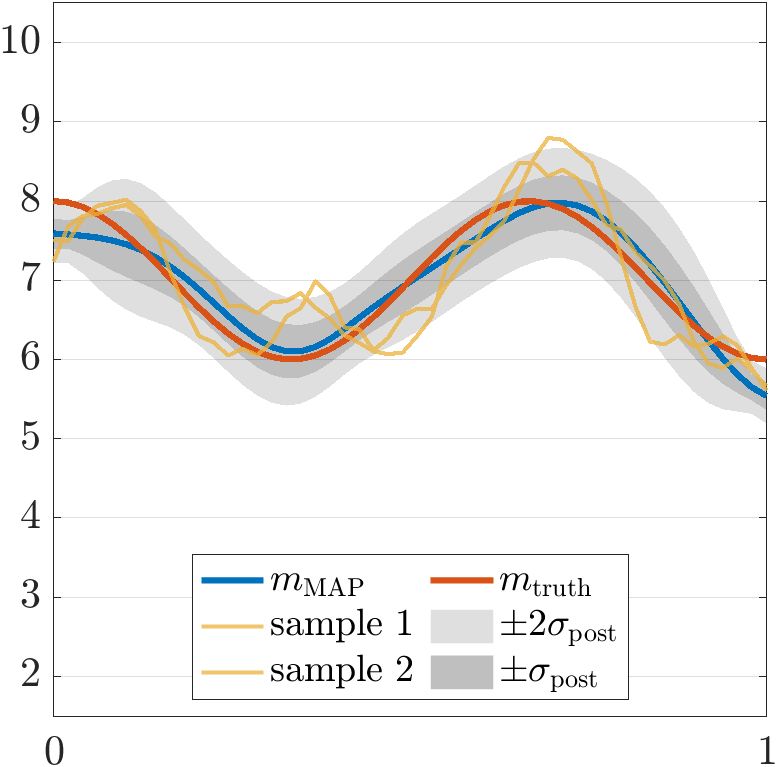}
\hfill\hfill\hfill
\includegraphics[height=0.27\textwidth]{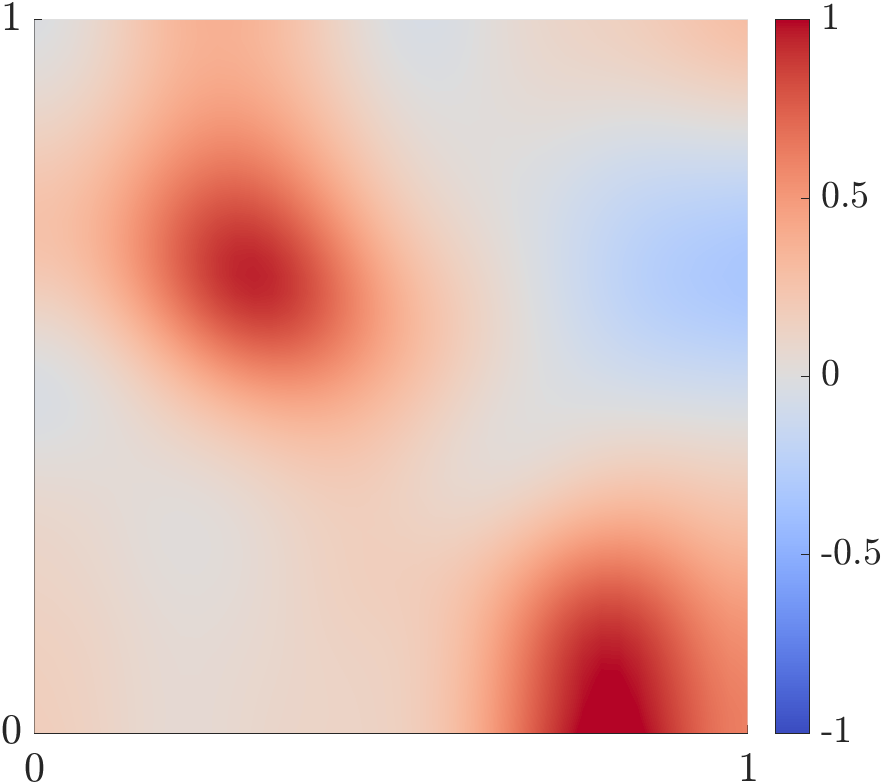}
\hfill
\includegraphics[height=0.27\textwidth]{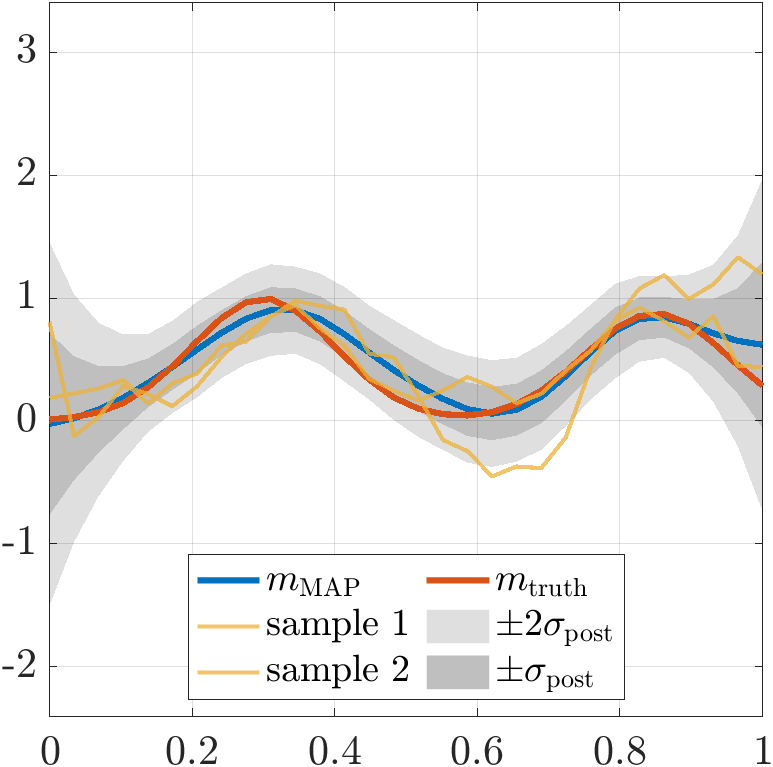}
    \caption{Posterior estimates found when using the true auxiliary parameter value during inversion. For Example 1 (left) shown are is the MAP estimate in blue, the truth in red, two samples from the (Gaussian approximation to the) posterior in yellow, while the plus/minus one and two standard deviation intervals,  $\pm\sigma$ and $\pm2\sigma$ are shaded in gray. For Example 2 we show the MAP estimate (centre) as well as the 
    one-dimensional marginal posterior plots (right) along the line from (0,1) to (1,0) (see Figure~\ref{fig: MAPQPAT}) with the MAP is shown in blue, the truth in shown in red, while two samples from the (Gaussian approximation to the) posterior are shown in yellow, and the plus/minus one and two standard deviation intervals,  $\pm\sigma$ and $\pm2\sigma$ are shaded in gray.}
    \label{fig: REFex1}
\end{figure}

Here we show posterior estimates of $\parm$ for both examples computed when the true value of the auxiliary parameter $\parb$ is used in the respective surrogate models, i.e., when we use $\cfwd(\parm):=\afwd(\parm,\parb_{\rm true})$ (which, in the nonlinear diffusion case, makes the forward model nonlinear in the state and thus significantly more involved). In Figure~\ref{fig: REFex1} we show the (Gaussian approximation to the) reference posteriors for Examples 1 and 2. These reference solutions serve as a baseline for the 2-Wasserstein metric comparisons in Section~\ref{sec: postests}. 

As a further reference, in Figure~\ref{fig: ex1SF} we show the posterior estimates found using the sample-free approach (discussed in Section~\ref{sec:sampFree}).

\begin{figure}[t!]
    \centering
\includegraphics[height=0.3\textwidth]{SF1CV.png}   
    \caption{The sample-free posterior estimate, with the MAP estimate shown in blue, the truth in red, two samples from the (Gaussian approximation to the) posterior in yellow, and the plus/minus one and two standard deviation intervals,  $\pm\sigma$ and $\pm2\sigma$ shaded in gray.}
    \label{fig: ex1SF}
\end{figure}

\end{document}